\documentclass[english,12pt]{smfart}
\usepackage{amssymb}
\usepackage{amsfonts}
\usepackage{amscd}
\usepackage[all]{xypic}
\usepackage[english]{babel}

\swapnumbers

\def\ac{{\rm ac}}

\def\11{{\mathbf 1}}
\def\ordjac{{\rm ord jac}}

\let\got\mathfrak
\def\gP{{\got P}}
\def\LO{{\cL_{\cO}}}

\def\chara{{\rm char}}

\def\Var{{\rm Var}}

\def\Def{{\rm Def}}

\def\Var{{\rm Var}}
\def\Def{{\rm Def}}
\def\GDef{{\rm GDef}}
\def\RDef{{\rm RDef}}

\let\got\mathfrak

\def\AA{{\mathbf A}}

\def\FF{{\mathbf F}}

\def\LL{{\mathbf L}}

\def\NN{{\mathbf N}}

\def\QQ{{\mathbf Q}}
\def\RR{{\mathbf R}}

\def\ZZ{{\mathbf Z}}

\def\cA{{\mathcal A}}
\def\cB{{\mathcal B}}
\def\cC{{\mathcal C}}

\def\cF{{\mathcal F}}

\def\cL{{\mathcal L}}
\def\cM{{\mathcal M}}

\def\cO{{\mathcal O}}
\def\cP{{\mathcal P}}

\def\cU{{\mathcal U}}

\def\cX{{\mathcal X}}
\def\cY{{\mathcal Y}}

\mathchardef\alphag="7C0B \mathchardef\betag="7C0C
\mathchardef\gammag="7C0D \mathchardef\deltag="7C0E
\mathchardef\varepsilong="7C22 \mathchardef\varphig="7C27
\mathchardef\psig="7C20 \mathchardef\zetag="7C10
\mathchardef\epsilong="7C0F \mathchardef\rhog="7C1A
\mathchardef\taug="7C1C \mathchardef\upsilong="7C1D
\mathchardef\iotag="7C13 \mathchardef\thetag="7C12
\mathchardef\pig="7C19 \mathchardef\sigmag="7C1B
\mathchardef\etag="7C11 \mathchardef\omegag="7C21
\mathchardef\kappag="7C14 \mathchardef\lambdag="7C15
\mathchardef\mug="7C16 \mathchardef\xig="7C18
\mathchardef\chig="7C1F \mathchardef\nug="7C17
\mathchardef\varthetag="7C23 \mathchardef\varpig="7C24
\mathchardef\varrhog="7C25 \mathchardef\varsigmag="7C26
\mathchardef\Omegag="7C0A \mathchardef\Thetag="7C02
\mathchardef\Sigmag="7C06 \mathchardef\Deltag="7C01
\mathchardef\Phig="7C08 \mathchardef\Gammag="7C00
\mathchardef\Psig="7C09 \mathchardef\Lambdag="7C03
\mathchardef\Xig="7C04 \mathchardef\Pig="7C05
\mathchardef\Upsilong="7C07

\newtheorem{theorem}[subsection]{Theorem}
\newtheorem{lem}[subsection]{Lemma}

\newtheorem{prop}[subsection]{Proposition}

\newtheorem{conj}[subsection]{Conjecture}

\theoremstyle{definition}
\newtheorem{definition}[subsection]{Definition}

\newtheorem{def-prop}[subsubsection]{Definition-Proposition}

\theoremstyle{remark}
\newtheorem{remark}[subsection]{Remark}

\theoremstyle{plain}

\numberwithin{equation}{subsection}

\def\boxit#1#2{\setbox1=\hbox{\kern#1{#2}\kern#1}%
\dimen1=\ht1 \advance\dimen1 by #1 \dimen2=\dp1
\advance\dimen2 by #1
\setbox1=\hbox{\vrule height\dimen1 depth\dimen2\box1\vrule}%
\setbox1=\vbox{\hrule\box1\hrule}%
\advance\dimen1 by .4pt \ht1=\dimen1 \advance\dimen2 by
.4pt \dp1=\dimen2 \box1\relax}

\newcommand{\sur}[2]{\genfrac{}{}{0pt}{}{#1}{#2}}
\renewcommand{\theequation}{\thesubsection.\arabic{equation}}

\let\got\mathfrak
\def\AA{{\mathbf A}}

\def\FF{{\mathbf F}}

\def\LL{{\mathbf L}}

\def\NN{{\mathbf N}}

\def\QQ{{\mathbf Q}}
\def\RR{{\mathbf R}}

\def\ZZ{{\mathbf Z}}

\def\LPas{\cL_{\rm DP}}
\def\LPre{\cL_{\rm DP,P}}

\def\cA{{\mathcal A}}
\def\cB{{\mathcal B}}
\def\cC{{\mathcal C}}

\def\cF{{\mathcal F}}

\def\cL{{\mathcal L}}
\def\cM{{\mathcal M}}

\def\cO{{\mathcal O}}
\def\cP{{\mathcal P}}

\def\cU{{\mathcal U}}

\def\cX{{\mathcal X}}
\def\cY{{\mathcal Y}}

\mathchardef\alphag="7C0B \mathchardef\betag="7C0C
\mathchardef\gammag="7C0D \mathchardef\deltag="7C0E
\mathchardef\varepsilong="7C22 \mathchardef\varphig="7C27
\mathchardef\psig="7C20 \mathchardef\zetag="7C10
\mathchardef\epsilong="7C0F \mathchardef\rhog="7C1A
\mathchardef\taug="7C1C \mathchardef\upsilong="7C1D
\mathchardef\iotag="7C13 \mathchardef\thetag="7C12
\mathchardef\pig="7C19 \mathchardef\sigmag="7C1B
\mathchardef\etag="7C11 \mathchardef\omegag="7C21
\mathchardef\kappag="7C14 \mathchardef\lambdag="7C15
\mathchardef\mug="7C16 \mathchardef\xig="7C18
\mathchardef\chig="7C1F \mathchardef\nug="7C17
\mathchardef\varthetag="7C23 \mathchardef\varpig="7C24
\mathchardef\varrhog="7C25 \mathchardef\varsigmag="7C26
\mathchardef\Omegag="7C0A \mathchardef\Thetag="7C02
\mathchardef\Sigmag="7C06 \mathchardef\Deltag="7C01
\mathchardef\Phig="7C08 \mathchardef\Gammag="7C00
\mathchardef\Psig="7C09 \mathchardef\Lambdag="7C03
\mathchardef\Xig="7C04 \mathchardef\Pig="7C05
\mathchardef\Upsilong="7C07

\DeclareMathOperator*{\Spec}{Spec}
\def\ord{{\rm ord} \, }

\begin{document}

\title[Ax-Kochen-Er{\v s}ov Theorems
and motivic integration]{Ax-Kochen-Er{\v s}ov Theorems for
$P$-adic integrals and motivic integration}

\author{Raf Cluckers}
\address{Katholieke Universiteit Leuven,
Departement wiskunde, Celestijnenlaan 200B, B-3001 Leu\-ven,
Bel\-gium. Current address: \'Ecole Normale Sup\'erieure,
D\'epartement de ma\-th\'e\-ma\-ti\-ques et applications, 45 rue
d'Ulm, 75230 Paris Cedex 05, France}
\email{raf.cluckers@wis.kuleuven.ac.be}

\author{Fran\c cois Loeser}
\address{{\'E}cole Normale Sup{\'e}rieure,
D{\'e}partement de math{\'e}matiques et applications,
45 rue d'Ulm,
75230 Paris Cedex 05, France
(UMR 8553 du CNRS)}
\email{Francois.Loeser@ens.fr}

%\begin{abstract}We express the Lefschetz number of iterates of the
%monodromy of a function on a smooth complex algebraic variety
%in terms of the Euler characteristic of a space of truncated arcs.
%\end{abstract}
\maketitle

\section{Introduction}This paper is concerned with extending classical
results \`a la Ax-Kochen-Er{\v s}ov to  $p$-adic integrals in a
motivic framework. The first section is expository, starting
from Artin's conjecture and the classical work of Ax, Kochen, and
Er{\v s}ov and ending with recent work of Denef and Loeser
giving a motivic version of the results of Ax, Kochen, and
Er{\v s}ov. In
that section we have chosen to adopt a quite informal style, since
the reader will find precise technical statements of more general
results in later sections. We also explain the cell decomposition
Theorem of Denef-Pas and how it leads to a quick proof of the
results of Ax, Kochen, and Er{\v s}ov. In sections \ref{cmf},
\ref{gm} and \ref{glob},
we present our new, general construction of motivic
integration, in the framework of constructible motivic functions.
This has been announced in \cite{cr1} and \cite{cr2} and is
developed in the paper \cite{cl}. In the last two sections we
explain the relation to $p$-adic integration and we announce
general Ax-Kochen-Er{\v s}ov Theorems for integrals depending on
parameters.
We conclude
the paper by discussing briefly the relevance of our
results to the study of
orbital integrals and the Fundamental Lemma.

The present text is an expanded and   updated version of a  talk given
by the senior author at the Miami Winter School
``Geometric Methods in Algebra and Number Theory''.
We  would like to thank the organizers for providing
such a  nice and congenial opportunity for  presenting our work.

\section{From Ax-Kochen-Er{\v s}ov to motives}\label{ak}
\subsection{Artin's conjecture}Let
$i$ and $d$ be integers. A field $K$ is said to be $C_i (d)$ if
every homogeneous polynomial of degree $d$ with coefficients in
$K$ in $d^i + 1$ (effectively appearing) variables has a non
trivial zero in $K$. Note we could replace ``in $d^i + 1$
variables'' by ``in at least $d^i + 1$ variables'' in that
definition. When  the field $K$ is $C_i (d)$ for every $d$ we say
it is $C_i$. For instance for a field $K$ to be $C_0$ means to be
algebraically closed, and all finite fields are $C_1$, thanks to
the Chevalley-Warning Theorem. Also, one can prove without much
trouble that if the field $K$ is $C_i$ then the fields $K (X) $
and $K ((X))$ are $C_{i + 1}$. It follows in particular that the
fields $\FF_q ((X))$ are $C_2$.

\begin{conj}[Artin]The $p$-adic fields $\QQ_p$ are
$C_2$.
\end{conj}

In 1965 Terjanian \cite{terj} gave an example of homogeneous form
of degree $4$ in $\QQ_2$ in $18 > 4^2$ variables having only
trivial zeroes in $\QQ_2$, thus giving a counterexample to Artin's
Conjecture. Let us briefly recall Terjanian's construction,
refering to \cite{terj} and \cite{delon} for more details. The
basic idea is the following: if $f$ is a homogeneous polynomial of
degree $4$ in $9$ variables with coeffiecients in $\ZZ$, such
that, for every $x$ in $\ZZ^9$, if $f (x) \equiv 0 \mod 4$, then
$2$ divides $x$, then the polynomial in $18$ variables $h (x, y) =
f (x) + 4 f (y)$ will have no non trivial zero in $\QQ_2$. An
example of such a polynomial $f$ is given  by
\begin{equation}f = n (x_1, x_2, x_3) +
n (x_4, x_5, x_6) +
n (x_7, x_8, x_9)
\end{equation}
with
\begin{equation}
n (X, Y, Z) =
X^2YZ +
XY^2Z +
XYZ^2
+X^2
+ Y^2
+ Z^2
- X^4 - Y^4 - Z^4.
\end{equation}

At about the same time, Ax and Kochen proved that, if not true,
Artin's conjecture is asymptotically true in the following sense:

\begin{theorem}[Ax-Kochen]\label{AK1}
An integer $d$ being fixed, all but finitely many fields
$\QQ_p$ are $C_2 (d)$.
\end{theorem}

\subsection{Some Model Theory}In fact, Theorem \ref{AK1} is a special instance of  the following,
much more general, statement:

\begin{theorem}[Ax-Kochen-Er{\v s}ov]\label{AK2}
Let $\varphi$ be a sentence in the language of rings. For all but
finitely prime numbers $p$, $\varphi$ is true in $\FF_p ((X))$ if
and only if it is true in $\QQ_p$. Moreover, there exists an
integer $N$ such that for any two local fields $K,K'$ with
isomorphic residue fields of characteristic $>N$ one has that
$\varphi$ is true in $K$ if and only if it is true in $K'$.
\end{theorem}

By a sentence in the language of rings, we mean a formula, without
free variables, built from symbols $0$, $+$, $-$, $1$, $\times$,
symbols for  variables, logical connectives $\wedge$, $\vee$,
$\neg$, quantifiers $\exists$, $\forall$ and the equality symbol
$=$. It is very important that in this language, any given natural
number can be expressed - for instance $3 $ as $1 + 1 + 1$-
 but that quantifiers
running for instance over natural numbers are not allowed. Given a
field $k$, we may interpret any such formula $\varphi$  in $k$, by
letting the quantifiers run over $k$, and, when $\varphi$ is a
sentence, we may say whether $\varphi$ is true in $k$ or not.
Since for a field to be $C_2(d)$ for a fixed $d$ may be expressed
by a sentence in the language of rings, we see that Theorem
\ref{AK1} is a special case of Theorem \ref{AK2}. On the other
hand, it is for instance impossible to express by a single
sentence in the language of rings for a field to be algebraically
closed.

In fact, it is natural to introduce here the language of valued
fields. It is a language with two sorts of variables. The first
sort of variables will run over the valued field  and the second
sort of variables will run over the value group. We shall use the
language of rings over the  valued field variables and the
language of ordered abelian groups $0, + , -, \geq$ over the value
group variables. Furthemore, there will be an additional
functional symbol $\ord$, going from the valued field sort  to the
value group sort, which will be interpreted as assigning to a non
zero element in the valued field its valuation.

\begin{theorem}[Ax-Kochen-Er{\v s}ov]\label{AKE}
Let $K$ and $K'$ be two henselian valued fields of residual
characteristic zero. Assume their residue fields $k$ and $k'$ and
their value groups $\Gamma$ and $\Gamma'$ are elementary
equivalent, that is, they have the same set of true sentences in
the rings, resp. ordered abelian groups, language. Then $K$ and
$K'$ are elementary equivalent, that is, they satisfy the same set
of formulas in the valued fields language.
\end{theorem}

We shall explain  a proof of Theorem 
\ref{AK2} after Theorem \ref{quantel}. Let us sketch how Theorem
\ref{AK2} also follows from Theorem \ref{AKE}. Indeed this follows
directly from the classical ultraproduct construction. Let
$\varphi$ be a given sentence in the language of valued fields.
Suppose by contradiction that for each $r$ in $ \NN$ there exist two
local fields $K_r,K_r'$ with isomorphic residue field of
characteristic $>r$ and such that $\varphi$ is true in $K_r$ and
false in $K'_r$. Let $U$ be a non principal ultrafilter on $\NN$.
Denote by $\FF_U$ the corresponding ultraproduct of the residue
fields of $K_r$, $r$ in $\NN$. It is a field of characteristic zero.
Now let $K_U$ and $K_U'$ be respectively the ultraproduct relative
to $U$ of the fields $K_r$ and $K'_r$. They are both henselian
with residue field $\FF_U$ and value group $\ZZ_U$, the
ultraproduct over $U$ of the ordered group $\ZZ$. Hence certainly
Theorem \ref{AKE} applies to $K_U$ and $K'_U$. By the very
ultraproduct construction, $\varphi$ is true in $K_U$ and false in
$K'_U$, which is a contradiction.

\subsection{}\label{dpcd}
In this paper, we shall in fact consider, instead of the language
of valued fields, what we call a language of Denef-Pas, $\LPas$.
It it is a language with 3 sorts, running respectively over valued
field, residue field, and value group variables. For the first 2
sorts, the language is the ring language and for the last sort, we
take any extension of the language of ordered abelian groups. For
instance, one may choose for the last sort the Presburger language
$ \{+, 0, 1, \leq\} \cup \{\equiv_n\ \mid n\in \NN,\ n
>1\},
$
where $\equiv_n$ denote equivalence modulo $n$. We denote the corresponding
Denef-Pas language by $\LPre$.
We also have two additional symbols, $\ord$ as before,
and a functional symbol $\ac$, going from the valued field sort to the residue field sort.

A typical example of a structure for that language is the field of
Laurent series $k((t))$ with the standard valuation $\ord : k
((t))^{\times} \rightarrow \ZZ$ and $\ac$ defined by $\ac (x) = x
t^{-\ord (x)} \mod t$ if $x \not=0$ in $k ((t))$ and by $\ac (0) =
0$.\footnote{Technically speaking, any function symbol of a first
order language must have as domain a product of sorts; a concerned
reader may choose an arbitrary extension of $\ord$ to the whole
field $K$; sometimes we will use ${\rm ord}_0:K\to \ZZ$ which
sends $0$ to $0$ and nonzero $x$ to $\ord(x)$.} Also, we shall
usually add to the language constant symbols in the first, resp.
second, sort for every element of $k ((t))$ resp. $k$, thus
considering formulas with coefficients in $k ((t))$, resp. $k$, in
the valued field, resp. residue field, sort. Similarly, any finite
extension of $\QQ_p$ is naturally a structure for that language,
once a uniformizing parameter $\varpi$ has been chosen; one just
sets $\ac (x) = x \varpi^{-\ord (x)} \mod \varpi$ and $\ac (0) =
0$. In the rest of the paper, for $\QQ_p$ itself, we shall always
take $\varpi = p$.

We now consider a valued field $K$ with residue field $k$ and
value group $\ZZ$. We assume $k$ is of characteristic zero, $K$ is
henselian and admits an angular component map, that is, a map $\ac
: K \rightarrow k$ such that $\ac (0) = 0$, $\ac$ restricts to a
multiplicative morphism $K^{\times} \rightarrow k^{\times}$, and
on the set $\{x \in K, \ord (x) = 0\}$, $\ac$ restricts to the
canonical projection to $k$. We also assume that $(K, k
, \Gamma, \ord, \ac)$ is a structure for the language $\LPas$.

We call a subset $C$ of $K^m \times k^n \times \ZZ^r$ definable if
it may be defined by a $\LPas$-formula. We call a function $h : C
\rightarrow K$ definable if its graph is definable.

\begin{definition} Let $D\subset K^m \times k^{n+1}\times\ZZ$ and $c:K^m \times k^n \to K$ be
definable. For $\xi$ in $k^n$, we set
\begin{equation*}
\begin{split}
A (\xi) =
\Big\{(x, t) \in K^m \times K & \Bigm\vert (x, \xi,\ac(t - c (x, \xi)),{\rm ord}_0(t - c (x, \xi))) \in D\},\\
\end{split}
\end{equation*}
where ${\rm ord}_0(x)=\ord(x)$ for $x\not=0$ and ${\rm
ord}_0(0)=0$. If for every $\xi$ and $\xi'$ in $k^n$ with $\xi
\not= \xi'$, we have $A (\xi) \cap A (\xi') = \emptyset$, then we
call
\begin{equation}
A = \bigcup_{\xi \in k^n} A (\xi)
\end{equation}
a cell in $K^m \times K$ with parameters $\xi$ and center $c (x, \xi)$.
\end{definition}

Now can state the following version of the cell decomposition Theorem of Denef and Pas:

\begin{theorem}[Denef-Pas \cite{Pas}]\label{cellth}Consider functions
$f_1 (x, t)$, \dots, $f_r (x, t)$ on $K^m \times K$ which are
polynomials in $t$ with coefficients definable functions from
$K^m$ to $K$. Then, $K^m \times K$ admits a finite partition into
cells $A$ with parameters $\xi$ and center $c (x, \xi)$, such
that, for every $\xi$ in $k^n$, $(x, t)$ in $A (\xi)$, and $1 \leq
i \leq r$, we have,
\begin{equation}
{\rm ord}_0 f_i (x, t) = {\rm ord}_0 h_i (x, \xi) (t - c (x,
\xi))^{\nu_i}
\end{equation}
and
\begin{equation}
\ac f_i (x, t) = \xi_{i},
\end{equation}
where the functions $h_i (x, \xi)$ are definable and $\nu_i,n$ are
in $\NN$ and where ${\rm ord}_0(x)=\ord(x)$ for $x\not=0$ and
${\rm ord}_0(0)=0$.
\end{theorem}

Using Theorem \ref{cellth} it is not difficult to prove by
induction on the number of valued field variables the following
quantifier elimination result (in fact, Theorems \ref{cellth} and
\ref{quantel} have a joint proof in \cite{Pas}):

\begin{theorem}[Denef-Pas \cite{Pas}]\label{quantel}Let $K$ be a valued
field satisfying the above conditions. Then, every formula in
$\LPas$ is equivalent to a formula without quantifiers running
over the valued field variables.
\end{theorem}

Let us now explain why Theorem \ref{AK2} follows easily from
Theorem \ref{quantel}. Let $U$ be a non principal ultrafilter on
$\NN$. Let $K_r$ and $K_r'$ be local fields for every  $r$ in $\NN$,
such that the residue field of $K_r$ is isomorphic to the residue
field of $K'_r$ and has characteristic $>r$. We consider again the
fields $K_U$ and $K'_U$ that are respectively the ultraproduct
relative to $U$ of the fields $K_r$ and $K_r'$. By the argument we
already explained it is enough to prove that these two fields are
elementary equivalent. Clearly they have isomorphic residue fields
and isomorphic value groups (isomorphic as ordered groups).
Furthermore they both satisfy the hypotheses of Theorem
\ref{quantel}. Consider a sentence true for $K_U$. Since it is
equivalent to a sentence with quantifiers running only over the
residue field variables and the value group variables, it will
also be true for $K_U'$, and vice versa.

Note that the use of cell decomposition to prove
Ax-Kochen-Er{\v s}ov type results goes back to
P.J. Cohen \cite{cohen}.

\subsection{From sentences to formulas}Let
$\varphi$ be a formula in the language of valued
fields or,  more generally, in the language $\LPre$ of Denef-Pas.
We assume that $\varphi$ has $m$ free
valued field variables  and no
free residue field nor value group variables.
For every valued field $K$ which is a structure for the language
$\LPas$,
we denote by $h_{\varphi} (K)$
the set of points $(x_1, \dots, x_m)$ in $K^m$
such that $\varphi (x_1, \dots, x_m)$ is true.

When $m = 0$, $\varphi$ is  a sentence and
$h_{\varphi} (K)$ is either the one point set or the empy set, depending
on whether $\varphi$ is true in $K$ or not.
Having  Theorem \ref{AK2} in mind, a natural question is to compare
$h_{\varphi} (\QQ_p)$ with $h_{\varphi} (\FF_p ((t)))$.

An answer is provided by the following statement:

\begin{theorem}[Denef-Loeser \cite{JAMS}]\label{thdl}Let
$\varphi$ be a formula in the language $\LPre$ with  $m$ free
valued field variables  and no free residue field nor value group
variables. There exists a virtual motive $M_{\varphi}$,
canonically attached to $\varphi$, such that, for almost all prime
numbers $p$, the volume of $h_{\varphi} (\QQ_p)$ is finite if and
only if the volume of $h_{\varphi} (\FF_p ((t)))$ is finite, and
in this case they are both equal to the number of points of
$M_{\varphi}$ in $\FF_p$.
\end{theorem}

Here we have chosen to state Theorem \ref{thdl} in an informal,
non technical way. A detailed presentation of more general results
we recently obtained is given in \S \kern .15em \ref{18}. A few
remarks are necessary in order to explain the statement of Theorem
\ref{thdl}. Firstly, what is meant by volume? Let $d$ be an
integer such that for almost all $p$, $h_{\varphi} (\QQ_p)$ is
contained in $X (\QQ_p)$, for some subvariety of dimension $d$ of
$\AA^m_{\QQ}$. Then the volume is taken with respect to the
canonical $d$-dimensional measure (cf. \S \kern .15em \ref{padic}
and \ref{18}). Implicit in the statement of the Theorem is the
fact that $h_{\varphi} (\QQ_p)$ and  $h_{\varphi} (\FF_p ((t)))$
are measurable (at least for almost all $p$ for the later one).
Originally, cf. \cite{JAMS} \cite{pek} \cite{dw}, the  virtual
motive $M_{\varphi}$ lies in a certain completion of the ring
$K_0^{\rm mot} ({\rm Var}_k) \otimes \QQ$ explained in \ref{4545}
(in particular, $K_0^{\rm mot} ({\rm Var}_k)$ is a subring of
the Grothendieck ring of Chow motives with rational coefficients),
but it now follows from the new construction of motivic
integration developed in \cite{cl}  that we can take $M_{\varphi}$
in the ring obtained from $K_0^{\rm mot} ({\rm Var}_k) \otimes
\QQ$ by inverting the Lefschetz motive $\LL$ and $1 - \LL^{- n}$
for $n > 0$.

One should note that even for $m = 0$, Theorem
\ref{thdl} gives more information than Theorem \ref{AK2}, since it
says that for almost all $p$  the validity of $\varphi$ in $\QQ_p$
and $\FF_p ((t))$ is governed by the virtual motive $M_{\varphi}$.
Finally, let us note that Theorem \ref{AK2} naturally extends to
integrals of definable functions as will be explained in \S \kern
.15em\ref{18}.

The proof of Theorem  \ref{thdl} is based on motivic integration.
In the next sections we shall give a quick overview of the  new
general construction of motivic integration given in \cite{cl},
that allows one to integrate a very general class of functions,
constructible motivic functions. These results have already been
announced in a condensed way in the notes \cite{cr1} and
\cite{cr2}; here, we are given the opportunity to present them
more leisurely and with some more details.

\section{Constructible motivic functions}\label{cmf}

\subsection{Definable subassignments}
Let $\varphi$ be a formula in the language $\LPre$ with
coefficients in $k ((t))$, resp. $k$, in the valued field, resp.
residue field, sort, having say respectively $m$, $n$, and $r$
free variables in the various sorts. To such  a formula $\varphi$
we assign, for every field $K$ containing $k$, the subset
$h_{\varphi} (K)$ of $K ((t))^m \times K^n \times \ZZ^r$
consisting of all points satisfying $\varphi$. We shall call the
datum of such subsets for all $K$ definable (sub)assignments. In
analogy with algebraic geometry, where the emphasis is not put
anymore on equations but on the functors they define, we consider
instead of formulas the corresponding subassignments (note $K
\mapsto h_{\varphi} (K)$ is in general not a functor).
Let us make these definitions more precise.

First, we recall the definition of subassignments, introduced in
\cite{JAMS}. Let $F : \cC \rightarrow {\rm Ens}$ be a functor from
a category  $\cC$ to the category of sets. By a subassignment $h$
of  $F$ we mean the datum, for every  object $C$ of $\cC$, of a
subset $h (C)$ of $F (C)$. Most of the standard operations of
elementary set theory extend trivially to subassignments. For
instance, given subassignments $h$ and $h'$ of the same functor,
one defines subassignments $h \cup h'$, $h \cap h'$ and the
relation $h \subset h'$, etc. When  $h \subset h'$ we say  $h$ is
a subassignment of $h'$. A morphism $f : h \rightarrow h'$ between
subsassignments of functors $F_1$ and $F_2$ consists of the datum
for every object $C$ of a map $f (C) : h (C) \rightarrow h'(C)$.
The graph of $f$ is the subassignment $C \mapsto {\rm graph} (f
(C))$ of $F_1 \times F_2$.

Next, we explain the notion of definable subassignments. Let
$k$ be a field and consider the category $F_k$ of fields
containing $k$. We denote by $h [m, n, r]$ the functor $F_k
\rightarrow {\rm Ens}$ given by  $h [m, n, r] (K) = K ((t))^m
\times K^n \times \ZZ^r$. In particular, $h [0, 0, 0]$ assigns the
one point set to every $K$.
%modif: We often...
%We often write $\ZZ^r$ for $h[0,0,r]$. 
To any formula $\varphi$ in
$\LPre$ with coefficients in $k ((t))$, resp. $k$, in the valued
field, resp. residue field, sort, having respectively $m$, $n$,
and $r$ free variables in the various sorts, we assign a
subsassignment $h_{\varphi}$ of $h [m, n, r]$, which associates to
$K$ in $F_k$ the subset $h_{\varphi} (K)$ of $h [m, n, r] (K)$
consisting of all points satisfying $\varphi$. We call such
subassignments definable subassignements. We denote by $\Def_k$
the   category whose objects are definable subassignments of some
$h [m, n, r]$, morphisms in $\Def_k$ being morphisms of
subassignments $f : h \rightarrow h'$ with  $h$ and $h'$ definable
subassignments of $h [m, n, r]$ and  $h [m', n', r']$ respectively
such that the graph of $f$ is a definable subassignment. Note that
$h [0, 0, 0]$ is the final object in this category.

If $S$ is an object of  $\Def_k$, we denote by  $\Def_S$ the
category of  morphisms $X \rightarrow S$ in $\Def_k$. If $f : X
\rightarrow S$ and $g : Y \rightarrow S$ are in $\Def_S$, we write
$X \times_S Y$ for the product in $\Def_S$ defined as $K \mapsto
\{(x, y) \in X (K) \times Y (K) \vert f (x) = g (y)\}$, with the
natural morphism to $S$. When $S = h [0, 0, 0)]$, we write $X
\times Y$ for $X \times_S Y$. We write  $S [m, n, r]$ for $S
\times h  [m, n, r]$, hence, $S [m, n, r] (K) = S (K) \times
K((t))^m \times K^n \times \ZZ^r$. By a  point $x$ of  $S$ we mean
a pair  $(x_0, K)$ with $K$ in $F_k$ and  $x_0$ a point of  $S
(K)$. We denote by  $\vert S \vert$ the set of  points of  $S$.
For such $x$ we then set $k (x) = K$. Consider a morphism $f : X
\rightarrow S$, with $X$ and $S$  respectively definable
subassignments of $h [m, n, r]$ and $h [m', n', r']$. Let $\varphi
(x, s)$ be a formula defining the graph of $f$ in $h [m + m', n +
n', r+ r']$. Fix a point $(s_0, K)$ of $S$. The formula $\varphi
(x, s_0)$ defines a subassignment in $\Def_{K}$. In this way we
get for $s$ a point of $S$ a functor  ``fiber at $s$'' $i_s^* :
\Def_S \rightarrow \Def_{k (s)}$.

\subsection{Constructible motivic functions}

In this subsection we define, for $S$ in $\Def_k$, the  ring
$\cC(S)$ of constructible motivic functions on $S$. The main goal
of this construction is that, as we will see in section \ref{gm},
motivic integrals with parameters in $S$ are constructible motivic
functions on $S$. In fact, in the construction of a measure, as we
all know since studying Lebesgue integration, positive functions
often play a basic fundamental role. This the reason why we also
introduce the semiring $\cC_+(S)$ of positive\footnote{Or maybe
better, non negative.} constructible motivic functions. A
technical novelty occurs here: $\cC(S)$ is the ring associated to
the semiring $\cC_+(S)$, but the canonical morphism $\cC_+ (S)
\rightarrow \cC (S)$ has in general no reason to be injective.

Basically, $\cC_+(S)$ and $\cC(S)$ are built up from
two kinds of functions.
The first type consists of elements of a certain Grothendieck (semi)ring.
Recall that in ``classical'' motivic integration as
developed in \cite{arcs}, the Grothendieck ring $K_0 (\Var_k)$
of algebraic
varieties over $k$ plays a key role. In the present setting the analogue
of the  category of algebraic
varieties over $k$ is the category of
definable subassignments of $h [0, n, 0]$, for some $n$, when $S = h[0,0,0]$.
Hence, for a general
$S$ in
$\Def_k$, it is natural to
consider the subcategory $\RDef_S$ of $\Def_S$ whose
objects are definable subassignments $Z$ of  $S \times h [0, n,
0]$, for variable $n$, the morphism $Z \rightarrow S$ being
induced by the  projection on  $S$.
The Grothendieck semigroup
$SK_0 (\RDef_S)$
is the quotient of the free semigroup on isomorphism classes of
objects $[Z \rightarrow S]$ in $\RDef_S$  by relations $[\emptyset
\rightarrow S] = 0$ and $[(Y \cup Y') \rightarrow S] + [(Y \cap
Y') \rightarrow S] = [Y \rightarrow S] + [Y' \rightarrow S]$. We
also denote by $K_0 (\RDef_S)$ the corresponding abelian group.
Cartesian product induces a unique semiring structure
on $SK_0 (\RDef_S)$, resp. ring structure on $K_0 (\RDef_S)$.

There are some easy functorialities.
For
every morphism $f : S \rightarrow S'$, there is a natural pullback
by $f^* : SK_0
(\RDef_{S'}) \rightarrow SK_0 (\RDef_S)$ induced by
the fiber product. If $f : S \rightarrow S'$ is a morphism
in  $\RDef_{S'}$,
composition with  $f$ induces a morphism $f_! : SK_0 (\RDef_S)
\rightarrow SK_0 (\RDef_{S'})$. Similar constructions apply to $K_0$.
That one can view elements  of $SK_0 (\RDef_{S})$
as functions on $S$ (which we even would like to integrate), is
illustrated in section \ref{padic} on $p$-adic integration and in
the introduction of \cite{cl}, in the part on integration
against Euler characteristic over the reals.

The second type of functions are certain
functions with values in the ring
\begin{equation}
A = \ZZ \Bigl[\LL, \LL^{-1}, \Bigl(\frac{1}{1 - \LL^{-i}}\Bigr)_{i
> 0}\Bigr],
\end{equation}
where, for the moment, $\LL$ is just considered as a symbol.
%modif: Note that...   
Note that a definable morphism $\alpha:S\to    h[0,0,1]$
determines a function $\vert S \vert \rightarrow \ZZ$, also
written $\alpha$, and a function $\vert S \vert \rightarrow A$
sending $x$ to $\LL^{\alpha(x)}$, written $\LL^\alpha$.
We consider the subring $\cP (S)$ of the ring of functions $\vert
S \vert \rightarrow A$ generated by constants in $A$ and by all
functions $\alpha$ and $\LL^\alpha$ with $\alpha:S \rightarrow
\ZZ$ definable morphisms. Now we should define positive functions
with values in $A$. For every real number  $q > 1$, let us  denote
by $\vartheta_q : A \rightarrow \RR$ the morphism sending $\LL$ to
$q$. We consider the subsemigroup $A_+$ of  $A$ consisting of
elements $a$ such that $\vartheta_q (a) \geq 0$ for all  $q > 1$
and we define $\cP_+ (S)$ as the semiring of functions in  $\cP
(S)$ taking their values in $A_+$.

Now we explain how to put together these two type of functions.
For  $Y$ a definable subassignment of $S$, we denote by ${\bf
1}_Y$ the function in $\cP (S)$ taking the value $1$ on $Y$ and
$0$ outside $Y$. We consider  the subring $\cP^0 (S)$ of  $\cP (S)
$, resp. the subsemiring $\cP^0_+ (S)$ of  $\cP_+ (S) $, generated
by functions of the form ${\bf 1}_Y$ with $Y$ a definable
subassignment of $S$, and by the constant function $\LL - 1$.
We have canonical morphisms $\cP^0 (S) \rightarrow K_0 (\RDef_S)$
and $\cP^0_+ (S) \rightarrow SK_0 (\RDef_S)$ sending ${\bf 1}_Y$
to $[Y \rightarrow S]$ and $\LL - 1$ to  the class of $S \times
(h[0, 1, 0] \setminus \{0\})$ in $K_0 (\RDef_S)$ and in  $SK_0
(\RDef_S)$, respectively. To simplify notation we shall denote by
$\LL$ and $\LL - 1$ the class of $S[0,1,0]$ and  $S \times (h[0,
1, 0] \setminus \{0\})$ in $K_0 (\RDef_S)$ and in  $SK_0
(\RDef_S)$.

We may now define the semiring of positive
constructible functions as
\begin{equation}
\cC_+ (S) = SK_0 (\RDef_S) \otimes_{\cP^0_+ (S)} \cP_+ (S)\end{equation}
and the ring of constructible functions as
\begin{equation}\cC (S) = K_0 (\RDef_S) \otimes_{\cP^0 (S)} \cP (S).
\end{equation}

If $f : S \rightarrow S'$ is a morphism in $\Def_k$, one shows in
\cite{cl} that the morphism $f^*$ may naturally be extended to a
morphism
\begin{equation}f^* : \cC_+ (S')
\longrightarrow \cC_+ (S).
\end{equation} If, furthermore, $f$ is a morphism in
$\RDef_{S'}$, one shows that the morphism $f_!$ may naturally be
extended to
\begin{equation}\label{fc}f_!
: \cC_+ (S) \longrightarrow \cC_+(S').\end{equation}
Similar functorialities exist for  $\cC$.

\subsection{Constructible motivic ``Functions''}\label{dim}
In fact, we shall need to consider not only functions as we just
defined, but functions defined almost everywhere in a given
dimension, that we call Functions.
(Note the capital in Functions.)

We start by defining a good notion of dimension for  objects of
$\Def_k$. Heuristically, that dimension corresponds to counting
the dimension only in the valued field variables, without taking
in account the remaining variables. More precisely, to any
algebraic subvariety $Z$ of $\AA^m_{k ((t))}$ we assign the
definable subassignment $h_Z$ of $h [m, 0, 0]$ given by  $h_Z (K)
= Z (K ((t)))$. The Zariski closure of a subassignment $S$ of $h
[m, 0, 0]$ is the intersection $W$ of all algebraic subvarieties
$Z$ of $\AA^m_{k ((t))}$ such that $S \subset h_Z$. We define the
dimension of $S$ as $\dim S := \dim W$. In the general case, when
$S$ is a subassignment of $h [m, n, r]$, we define  $\dim S$ as
the  dimension of the image of  $S$ under the  projection $h [m,
n, r] \rightarrow h [m, 0, 0]$.

One can prove, using
Theorem \ref{cellth}
and  results of van den Dries \cite{vdDries},
the following result, which is by no means obvious:

\begin{prop}\label{pdim}Two isomorphic objects of $\Def_k$
have the same dimension.
\end{prop}

For every non negative integer $d$, we denote by
$\cC_+^{\leq d} (S)$ the ideal of  $\cC_+ (S)$ generated by
functions $\11_{Z}$ with  $Z$ definable subassignments of $S$ with
$\dim Z \leq d$. We set  $C_+ (S) = \oplus_d  C^d_+ (S)$ with
$C^d_+ (S) := \cC_+^{\leq d} (S) / \cC_+^{\leq d-1} (S)$. It is a
graded abelian semigroup, and also a $\cC_+ (S)$-semimodule.
Elements of $C_+ (S) $ are called positive constructible Functions
on  $S$. If $\varphi$ is a function lying in $\cC_+^{\leq d} (S)$
but not in  $\cC_+^{\leq d - 1} (S)$, we denote by  $[\varphi]$
its  image in $C^d_+ (S)$.
One defines
similarly $C (S)$ from $\cC (S)$.

One of the reasons why we consider functions which are defined
almost everywhere originates in the differentiation of functions
with respect to the valued field variables: one may show that a
definable function $c : S\subset h[m,n,r] \rightarrow h [1, 0, 0]$
is differentiable (in fact even analytic) outside a definable
subassignment of $S$ of dimension $< {\rm dim} S$.
In particular, if $f : S \rightarrow S'$ is an isomorphism in
$\Def_k$, one may define
a function $\ordjac f$, the order of the jacobian of $f$, which is defined
almost everywhere and is equal almost everywhere to a definable
function, so we may define $\LL^{- \ordjac f}$ in $C_+^d (S)$ when
$S$ is of dimension $d$.
In \ref{dft}, we 
shall define $\LL^{- \ordjac f}$ using differential forms.

\section{Construction of the general motivic measure}\label{gm}

Let $k$ be a field of characteristic zero. Given
$S$ in  $\Def_k$,
we define
$S$-integrable Functions and construct pushforward morphisms
for these:

\begin{theorem}\label{mt}Let $k$ be a field of characteristic zero and
let $S$ be in  $\Def_k$.
There exists a unique functor $Z \mapsto
{\rm I}_S C_+(Z)$ from  $\Def_S$ to the category of abelian semigroups,
the functor of  $S$-integrable Functions,
assigning to every morphism $f : Z \rightarrow Y$ in
$\Def_S$ a morphism  $f_! : {\rm I}_S C_+(Z) \rightarrow {\rm I}_S C_+(Y)$
such that
for every  $Z$ in $\Def_S$, ${\rm I}_S C_+(Z)$ is a
graded subsemigroup of  $C_+ (Z)$  and  ${\rm I}_S C_+(S) = C_+ (S)$,
satisfying the following
list of axioms (A1)-(A8).
\end{theorem}

\noindent {\bf {\rm (A1a)} (Naturality)}

\noindent If  $S \rightarrow S'$ is a morphism in  $\Def_k$
and $Z$ is an object in $\Def_S$, then
any $S'$-integrable Function $\varphi$ in $C_+(Z)$ is $S$-integrable
and  $f_! (\varphi)$ is the same,
considered
in  ${\rm I}_{S'}$ or in ${\rm I}_{S}$.

\bigskip

\noindent {\bf {\rm (A1b)} (Fubini)}

\noindent A positive Function  $\varphi$ on $Z$ is $S$-integrable if and only if it is
$Y$-integrable and $f_! (\varphi)$ is
$S$-integrable.

\bigskip

\noindent {\bf {\rm (A2)} (Disjoint union)}

\noindent If $Z$ is the disjoint union of two definable
subassignments $Z_1$ and  $Z_2$, then the isomorphism $C_+ (Z)
\simeq C_+ (Z_1) \oplus C_+ (Z_2)$ induces an  isomorphism ${\rm
I}_S C_+ (Z) \simeq {\rm I}_S C_+ (Z_1) \oplus {\rm I}_S C_+
(Z_2)$, under which  $f_! = f_{|Z_1 !} \oplus f_{|Z_2 !}$.

\bigskip

\noindent {\bf {\rm (A3)} (Projection formula)}

\noindent For every $\alpha$ in  $\cC_+ (Y)$ and every  $\beta$
in ${\rm I}_S C_+ (Z)$, $\alpha f_! (\beta)$ is $S$-integrable if
and only if $f^* (\alpha) \beta$ is, and then $f_! (  f^* (\alpha)
\beta)  = \alpha f_! (\beta)$.

\bigskip

\noindent {\bf {\rm (A4)} (Inclusions)}

\noindent If  $i : Z \hookrightarrow Z'$ is the inclusion of
definable subassignments of the same object of $\Def_S$, $i_!$ is
induced by extension by zero outside  $Z$ and sends injectively
${\rm I}_S C_+ (Z)$ to ${\rm I}_S C_+ (Z')$.

\bigskip

\noindent {\bf {\rm (A5)} (Integration along
residue field variables)}

\noindent Let $Y$ be an object of  $\Def_S$ and denote by $\pi$ the projection
$Y [0, n, 0] \rightarrow Y$. A Function
$[\varphi]$ in  $C_+ (Y [0, n, 0])$ is $S$-integrable if and only if,
with  notations of \ref{fc},
$[\pi_! (\varphi)]$ is $S$-integrable and then
$\pi_! ([\varphi]) = [\pi_! (\varphi)]$.

Basically this  axiom means that
integrating with respect to
variables in the residue field
just amounts
to taking the  pushforward induced by composition at the level
of Grothendieck semirings.

\bigskip

\noindent {\bf {\rm (A6)} (Integration along
$\ZZ$-variables)}
Basically, integration along $\ZZ$-variables corresponds to summing
over the integers, but to state precisely (A6), we need to
perform some preliminary  constructions.

Let us consider a function in $\varphi$ in
$\cP (S [0, 0, r])$, hence
$\varphi$ is a function $\vert S \vert \times \ZZ^r \rightarrow A$.
We shall say  $\varphi$ is
$S$-integrable if for every
$q > 1$ and every  $x$ in $\vert S\vert$,
the series $\sum_{i \in \ZZ^r} \vartheta_q (\varphi (x, i))$
is summable. One proves that
if $\varphi$ is
$S$-integrable there exists a unique function $\mu_S (\varphi)$
in $\cP (S)$ such that
$\vartheta_q (\mu_S (\varphi) (x))$ is equal to the sum of the previous series
for all $q > 1$ and all  $x$ in $\vert S\vert$.
We denote by  ${\rm I}_S \cP_+ (S[0, 0, r])$
the set of  $S$-integrable functions in
$ \cP_+ (S[0, 0, r])$
and we set
\begin{equation}
{\rm I}_S \cC_+ (S[0, 0, r]) = \cC_+ (S)  \otimes_{\cP_+ (S)} {\rm I}_S \cP_+ (S[0, 0, r]).
\end{equation}
Hence ${\rm I}_S \cP_+ (S[0, 0, r])$
is a sub-$\cC_+ (S)$-semimodule of
$ \cC_+ (S[0, 0, r])$ and $\mu_S$ may be extended by tensoring to
\begin{equation}
\mu_S : {\rm I}_S \cC_+ (S[0, 0, r]) \rightarrow \cC_+ (S).
\end{equation}

Now we can state (A6):

\noindent Let $Y$ be an object of $\Def_S$ and denote by $\pi$ the  projection
$Y [0, 0, r] \rightarrow Y$. A Function
$[\varphi]$ in  $C_+ (Y [0, 0, r])$ is  $S$-integrable if and only if
there exists  $\varphi'$ in
$\cC_+ (Y [0, 0, r])$
with  $[\varphi'] = [\varphi]$ which is  $Y$-integrable in the previous
sense and such that  $[\mu_Y (\varphi')]$ is $S$-integrable. We then have
$\pi_! ([\varphi]) = [\mu_Y (\varphi')]$.

\bigskip

\noindent {\bf {\rm (A7)} (Volume of balls)}
It is natural to require  (by analogy with the $p$-adic case)
that the volume of a ball
$\{ z \in h [1, 0, 0] \vert \ac (z -c ) = \alpha , \ac (z -c ) = \xi\}$,
with $\alpha$ in $\ZZ$, $c$ in $k ((t))$ and $\xi$ non zero in $k$,
should be $\LL^{- \alpha -1}$.
(A7) is a relative version of that statement:

\noindent Let  $Y$ be an object in $\Def_S$ and let  $Z$ be the
definable subassignment of $Y [1, 0, 0] $ defined by  $\ord (z - c
(y)) = \alpha (y)$ and $\ac (z - c (y)) = \xi (y)$, with  $z$ the
coordinate on the $\AA^1_{k ((t))}$-factor and  $\alpha, \xi, c$
definable functions on $Y$ with values respectively in $\ZZ$,
$h[0, 1, 0] \setminus \{0\}$, and  $h[1, 0, 0]$. We denote by  $f
: Z \rightarrow Y$ the morphism induced by  projection. Then
$[\11_Z]$ is $S$-integrable if and only if $\LL^{-\alpha - 1}
[\11_Y]$ is, and then $f_! ([\11_Z]) = \LL^{-\alpha - 1} [\11_Y]$.

\bigskip

\noindent {\bf {\rm (A8)} (Graphs)} This last axiom expresses the
pushforward for graph projections.
It relates volume and differentials and is a special
case of the change of variables Theorem \ref{cvf}.

\noindent Let $Y$ be in  $\Def_S$ and let $Z$ be the definable
subassignment of $Y [1, 0, 0] $ defined by  $z - c (y) = 0$ with
$z$ the coordinate on the $\AA^1_{k ((t))}$-factor and  $c$ a
morphism  $Y \rightarrow h[1, 0, 0]$. We denote by  $f : Z
\rightarrow Y$ the morphism induced by projection. Then $[\11_Z]$
is $S$-integrable if and only if $\LL^{(\ordjac f) \circ f^{-1}}$ is,
and then  $f_! ([\11_Z]) = \LL^{(\ordjac f) \circ f^{-1}}$.

\bigskip

Once
Theorem \ref{mt} is proved, one may proceed as follows
to extend the constructions
from $C_+$ to $C$ .
One
defines ${\rm I}_S C (Z)$ as the subgroup of $C (Z)$ generated
by the image of  ${\rm I}_S C_+ (Z)$. One shows that if  $f : Z
\rightarrow Y$ is a  morphism in $\Def_S$, the  morphism $f_! :
{\rm I}_S C_+(Z) \rightarrow {\rm I}_S C_+(Y)$ has a natural
extension $f_! : {\rm I}_S C(Z) \rightarrow {\rm I}_S C(Y)$.

The relation
of Theorem \ref{mt} with motivic integration is the following.
When  $S$ is equal to  $h [0, 0, 0]$, the final object of
$\Def_k$, one writes ${\rm I} C_+ (Z)$ for
 ${\rm I}_S C_+ (Z)$ and we shall say integrable for $S$-integrable,
and similarly for  $C$.
Note that  ${\rm I} C_+ (h [0, 0, 0]) = C_+ (h [0, 0, 0]) =
SK_0 (\RDef_k) \otimes_{\NN [\LL - 1]} A_+$ and that
${\rm I}C (h [0, 0, 0]) =  K_0 (\RDef_k) \otimes_{\ZZ [\LL]} A$.
For $\varphi$ in  ${\rm I} C_+ (Z)$, or in  ${\rm I} C (Z)$, one defines the
motivic integral
$\mu (\varphi)$ by $\mu (\varphi) = f_! (\varphi)$ with
$f$ the morphism  $Z \rightarrow h [0, 0, 0]$.
Working in the  more general framework of Theorem \ref{mt}
to construct $\mu$
appears to be very convenient for inductions
occuring in the proofs. Also, it is not clear how to characterize
$\mu$ alone by existence and unicity properties.
Note also, that  one reason for the statement of
Theorem \ref{mt} to look somewhat cumbersone, is that we have to define
at once the notion of integrability and the value of the
integral.

\medskip

The proof of Theorem \ref{mt} is quite long and involved. In a
nutshell, the basic idea is the following. Integration along
residue field variables is controlled by (A5) and integration
along $\ZZ$-variables by (A6). Integration along valued field
variables is constructed one variable after the other. To
integrate with respect to one valued field variable, one may,
using (a variant of) the cell decomposition Theorem \ref{cellth}
(at the cost of introducing additional new residue field and
$\ZZ$-variables), reduce to the case of cells which is covered by
(A7) and (A8).
An important step is to show that this is independent of the choice of a
cell decomposition. When one integrates with respect to more than
one valued field variable (one after the other) it is crucial to
show that it is independent of the order of the variables, for
which we use a notion of bicells.

\medskip

In this new framework, we have the  following general form of the
change of variables Theorem, generalizing the corresponding
statements in \cite{arcs} and \cite{JAMS}.

\begin{theorem}\label{cvf}Let $f : X \rightarrow Y$ be an isomorphism between definable
subassignments of dimension $d$.
For every function  $\varphi$
in $\cC^{\leq d}_+(Y)$ having a non zero class in
$C^d_+ (Y)$, $[f^* (\varphi)]$ is $Y$-integrable
and
$f_! [f^* (\varphi)] = \LL^{(\ordjac f )\circ f^{-1}} [\varphi]$.
A similar statement holds in $C$.
\end{theorem}

\subsection{Integrals depending on parameters}
One pleasant feature of Theorem \ref{mt} is that it generalizes
readily to the relative setting of integrals depending on
parameters.

Indeed, let us  fix
$\Lambda$ in $\Def_k$ playing the role of a parameter space.
For $S$ in $\Def_{\Lambda}$,
we consider the ideal $\cC^{\leq d} (S \rightarrow \Lambda)$
of  $\cC_+ (S)$ generated by
functions $\11_{Z}$ with  $Z$ definable subassignment of $S$
such that all fibers of $Z \rightarrow \Lambda$ are of  dimension $\leq d$.
We set
\begin{equation}C_+ (S \rightarrow \Lambda) =
\bigoplus_d  C^d_+ (S \rightarrow \Lambda)
\end{equation}
with
\begin{equation}
C^d_+ (S \rightarrow \Lambda) := \cC^{\leq d}_+ (S \rightarrow \Lambda)  / \cC_+^{\leq d-1} (S \rightarrow \Lambda).
\end{equation}
It is a graded abelian semigroup (and also a
$\cC_+ (S)$-semimodule).  If $\varphi$
belongs to  $\cC_+^{\leq d} (S \rightarrow \Lambda)$
but not to  $\cC_+^{\leq d - 1} (S \rightarrow \Lambda)$,
we write $[\varphi]$ for its image in $C_+^d (S \rightarrow \Lambda)$.
The following relative analogue of Theorem  \ref{mt} holds.

\begin{theorem}\label{mtr}Let  $k$ be a field of characteristic zero, let
$\Lambda$ be in $\Def_k$, and let  $S$ be in $\Def_{\Lambda}$.
There exists a unique functor $Z \mapsto {\rm I}_S C_+(Z
\rightarrow \Lambda)$ from  $\Def_S$ to the category of  abelian
semigroups, assigning to every morphism $f : Z \rightarrow Y$ in
$\Def_S$ a morphism  $f_{! \Lambda} : {\rm I}_S C_+(Z\rightarrow
\Lambda)) \rightarrow {\rm I}_S C_+(Y\rightarrow \Lambda))$
satisfying properties analogue to   (A0)-(A8) obtained by
replacing $ C_+ (\_)$ by $C_+ (\_ \rightarrow \Lambda)$ and $\ordjac$
by its relative analogue $\ordjac_{\Lambda}$\footnote{Defined similarly
as $\ordjac$, but using relative differential forms.}.
\end{theorem}

Note that $C_+ (\Lambda
\rightarrow \Lambda) = \cC_+ (\Lambda)$ (and also
${\rm I}_{\Lambda} C_+ (\Lambda
\rightarrow \Lambda) = C_+ (\Lambda
\rightarrow \Lambda)$.
Hence, given
$f : Z \rightarrow \Lambda$ in
$\Def_{\Lambda}$, we may define the relative motivic measure
with respect to $\Lambda$ as the morphism
\begin{equation}
\mu_{\Lambda} :=
f_{!\Lambda } : {\rm I}_{\Lambda} C_+ (Z \rightarrow \Lambda)
\longrightarrow  \cC_+ (\Lambda).
\end{equation}

By the following statement,
$\mu_{\Lambda}$ indeed corresponds to integration
along the fibers over $\Lambda$:

\begin{prop}\label{ll}Let $\varphi$
be a  Function in  $C_+(Z \rightarrow \Lambda)$.
It belongs to
${\rm I}_{\Lambda} C_+(Z \rightarrow \Lambda)$ if and only if
for every point  $\lambda$ in  $\Lambda$, the  restriction $\varphi_{\lambda}$
of
$\varphi$ to the fiber  of  $Z$ at $\lambda$
is integrable. The motivic integral of
$\varphi_{\lambda}$
is then equal to
$i_{\lambda}^* (\mu_{\Lambda} (\varphi))$, for every $\lambda$ in
$\Lambda$.
\end{prop}

Similarly as in the absolute case,
one can also define the relative analogue
$C (S \rightarrow \Lambda)$ of $C (S)$,
and extend the notion of integrability and the  construction of
$f_{!\Lambda}$ to this setting.

\section{Motivic integration in a global setting  and comparison
with previous constructions}\label{glob}

\subsection{Definable subassignments on varieties}Objects
of $\Def_k$ are by construction affine, being subassignments of
functors $h [m, n, r] : F_k \rightarrow {\rm Ens}$ given by
$K \mapsto K((t))^m \times K^n \times \ZZ^r$.
We shall now consider their global analogues and extend the previous
constructions to the global setting.

Let  $\cX$ be a variety over $k ((t))$, that is, a reduced and
separated scheme of finite type over $k ((t))$, and let  $X$ be a
variety over $k$. For $r$ an integer  $\geq 0$, we denote by  $h
[\cX, X, r]$ the functor $F_k \rightarrow {\rm Ens}$ given by $K
\mapsto \cX (K ((t))) \times X (K) \times \ZZ^r$. When  $X = \Spec
k$ and $r = 0$, we write $h [\cX]$ for $h [\cX, X, r]$. If  $\cX$
and $X$ are affine and if  $i : \cX \hookrightarrow \AA^m_{k
((t))}$ and $j : X \hookrightarrow \AA^n_k$ are closed immersions,
we say a subassignment $h$ of  $h [\cX, X, r]$ is definable if its
image by the  morphism $h [\cX, X, r] \rightarrow h [m, n, r]$
induced by $i$ and $j$ is a definable subassignment of $h [m, n,
r]$. This definition does not depend on  $i$ and  $j$. More
generally, we shall say a subassignment $h$ of $h [\cX, X, r]$ is
definable if there exist coverings $(\cU_i)$ and  $(U_j)$ of
$\cX$ and $X$ by affine open subsets such that $h \cap h [\cU_i,
U_j, r]$ is a definable subassignment of  $h [\cU_i, U_j, r]$ for
every $i$ and $j$. We get in this way a category $\GDef_k$ whose
objects are definable subassignments of some $h [\cX, X, r]$,
morphisms being definable morphisms, that is, morphisms whose
graphs are definable subassignments.

The category $\Def_k$ is a full subcategory of $\GDef_k$.
Dimension as defined in \ref{dim} may be directly generalized to
objects of  $\GDef_k$ and Proposition \ref{pdim} still holds in
$\GDef_k$. Also, if $S$ is an object in $\GDef_k$, our definitions
of $\RDef_S$, $\cC_+ (S)$,
$\cC (S)$, $C_+ (S)$  and $C (S)$ extend.

\subsection{Definable differential forms and volume forms}\label{dft}In
the global setting, one does not integrate functions anymore, but
volume forms. Let us start by introducing differential forms in
the definable framework. Let  $h$ be a definable subassignment of
some $h [\cX, X, r] $. We denote by $\cA (h)$ the ring of
definable morphisms $h  \rightarrow h [\AA^1_{k ((t))}]$. Let us
define, for $i$ in $\NN$, the  $\cA (h)$-module $\Omega^i (h)$ of
definable $i$-forms on  $h$. Let $\cY$ be the closed subset of
$\cX$, which is the Zariski closure of the image of  $h$ under the
projection $\pi : h  [\cX, X, r] \rightarrow h [\cX] $. We denote
by $\Omega^i_{\cY}$ the sheaf of algebraic  $i$-forms on $\cY$, by
$\cA_{\cY}$ the  Zariski sheaf associated to the presheaf $U
\mapsto \cA (h [U])$ on $\cY$,  and by $\Omega^i_{h [\cY]}$ the
sheaf $\cA_{\cY} \otimes_{\cO_{\cY}} \Omega^i_{\cY}$. We set
\begin{equation}\Omega^i (h) :=
\cA (h) \otimes_{\cA (h [\cY])}    \Omega^i_{h [\cY]} (\cY),
\end{equation} the
$\cA (h [\cY])$-algebra structure on $\cA (h)$ given by
composition with  $\pi$.

We now assume $h$ is of  dimension $d$. We denote by  $\cA^< (h)$
the ideal of functions in  $\cA (h)$ that are zero outside a
definable subassignment of dimension $<d$.
There is a canonical morphism of abelian semi-groups $\lambda :
\cA (h) / \cA^< (h) \rightarrow C^d_+ (h)$ sending the class of a
function  $f$ to the class of $\LL^{- \ord f}$, with the
convention $\LL^{- \ord 0} = 0$.
We set $\tilde \Omega^d (h) = \cA (h) / \cA^< (h) \otimes_{\cA
(h)} \Omega^d (h) $, and we define the set $\vert \tilde \Omega
\vert_+ (h) $ of definable positive volume forms as the  quotient
of the free abelian semigroup on symbols $(\omega, g)$ with
$\omega$ in $\tilde \Omega^d (h) $ and $g$ in $C^d_+ (h)$ by
relations $(f \omega, g) = (\omega, \lambda (f) g)$, $(\omega, g +
g') = (\omega, g) + (\omega, g')$ and $(\omega, 0)= 0$, for $f$ in
$\cA (h) / \cA^< (h)$. We write $g |\omega|$ for the class
$(\omega, g)$, in order to have $g \vert f \omega \vert = g \LL^{-
\ord f} \vert \omega \vert$.
The $\cC_+ (h)$-semimodule structure on $C^d_+ (h)$ induces  after
passing to the quotient a  structure of semiring on $C^d_+ (h)$
and $\vert \tilde \Omega  \vert_+(h)$ is naturally endowed with a
structure of $C^d_+ (h)$-semimodule. We shall call an element
$\vert \omega \vert$ in $\vert \tilde \Omega  \vert_+ (h)$ a gauge
form if it is a generator of that semimodule. One should note
that in the present setting
gauge forms  always exist,
which is certainly not the case in the usual framework of
algebraic geometry.
Indeed,
gauge forms always exist locally (that is, in suitable affine charts),
and in our definable world there is no difficulty
in gluing local gauge forms to global ones.
One may define similarly $\vert \tilde
\Omega \vert(h)$, replacing  $C^{d}_+$ by $C^{d}$, but we shall
only consider $\vert \tilde \Omega  \vert_+(h)$ here.

If $h$ is definable subassignment of  dimension $d$ of $h [m, n,
r]$, one may construct, similarly as Serre \cite{serre} in the
$p$-adic case, a canonical gauge form $\vert \omega_0\vert_h$ on
$h$. Let us denote by  $x_1$, \dots, $x_m$ the coordinates on
$\AA^m_{k((t))}$  and consider the   $d$-forms $\omega_I :=
dx_{i_1} \wedge \dots \wedge dx_{i_d}$ for $I = \{i_1, \dots,
i_d\} \subset \{1, \dots m\}$, $i_1 < \dots < i_d$, and their
image $\vert \omega_I \vert_h$ in $\vert \tilde \Omega \vert_+
(h)$.
One may check there exists a unique element  $\vert
\omega_0\vert_h$ of $\vert \tilde \Omega \vert_+ (h)$, such that,
for every $I$, there exists definable functions with integral
values $\alpha_I$, $\beta_I$ on $h$, with $\beta_I$ only taking as
values $1$ and $0$, such that $\alpha_I+\beta_I>0$ on $h$, $\vert
\omega_I \vert_h = \beta_I \LL^{- \alpha_I} \vert \omega_0\vert_h$
in $\vert \tilde \Omega \vert_+ (h)$, and such that $\inf_I
\alpha_I = 0$.
%modif outside a definable subassignment of dimension $<d$.

If  $f : h \rightarrow h'$ is a morphism in $\GDef_k$ with $h$ and
$h'$ of dimension $d$ and all fibers of dimension $0$, there is a
mapping $f^* : \vert \tilde \Omega \vert_+ (h')  \rightarrow \vert
\tilde \Omega \vert_+ (h)$ induced by pull-back of differential
forms. This follows from the fact that $f$ is ``analytic'' outside
a definable subassignment of dimension $d - 1$ of $h$. If,
furthermore, $h$ and  $h'$ are  objects in $\Def_k$, one
defines  $\LL^{ -\ordjac f}$ by
\begin{equation}
f^* \vert
\omega_0\vert_{h'} = \LL^{ -\ordjac f} \vert \omega_0\vert_h.
\end{equation}

If  $\cX$ is a  $k((t))$-variety of  dimension $d$, and  $\cX^0$
is a $k [[t]]$-model of $\cX$, it is possible to define an element
$|\omega_0|$ in $|\tilde \Omega|_+ (h [\cX])$, which depends only
on $\cX^0$, and which is characterized by the following property:
for every open $U^0$ of  $\cX^0$ on which the  $k [[t]]$-module
$\Omega^d_{U^0 | k [[t]]} (U^0)$ is generated by a nonzero form
$\omega$, $|\omega_0|_{|h [U^0 \otimes \Spec k((t))]} = |\omega|$
in $|\tilde \Omega|_+(h [U^0 \otimes \Spec k((t))])$.

\subsection{Integration of volume forms and Fubini Theorem}
Now we are ready to construct motivic integration for volume
forms. In the affine case, using canonical gauge forms, one may
pass from volume forms to Functions in top dimension, and vice
versa. More precisely, let $f : S \rightarrow S'$ be a  morphism
in  $\Def_k$, with  $S$ of dimension $s$ and  $S'$ of dimension
$s'$. Every positive form $\alpha$ in $|\tilde \Omega|_+ (S)$ may
be written  $\alpha = \psi_{\alpha} |\omega_0|_S$ with
$\psi_{\alpha}$ in $C^s_+ (S)$. We shall say $\alpha$ is
$f$-integrable if $\psi_{\alpha}$ is $f$-integrable and we  then
set
\begin{equation}f_!^{\rm top} (\alpha) := \{f_! (\psi_{\alpha})\}_{s'}|\omega_0|_{S'},
\end{equation}
$\{f_! (\psi_{\alpha})\}_{s'}$ denoting the component of
$f_! (\psi_{\alpha})$ lying in  $C^{s'}_+ (S')$.

Consider now a morphism
$f : S \rightarrow S'$  in $\GDef_k$. The previous construction may be
globalized as follows.
Assume there exist  isomorphisms
$\varphi : T \rightarrow S$ and  $\varphi' : T' \rightarrow S'$
with  $T$ and  $T'$ in  $\Def_k$.
We denote by $\tilde f$ the morphism $ T \rightarrow T'$
such that $\varphi' \circ \tilde f = f \circ \varphi$.
We shall say  $\alpha$ in
$|\tilde \Omega|_+ (S)$ is $f$-integrable if
$\varphi^* (\alpha)$ is $\tilde f$-integrable
and we define then  $f^{\rm top}_! (\alpha)$
by the relation
\begin{equation}
{\tilde f}^{\rm top}_! (\varphi^* (\alpha))
= \varphi'{}^* (f^{\rm top}_! (\alpha)).
\end{equation}
It follows from Theorem \ref{cvf} that this definition is
independent of  the choice of the isomorphisms $\varphi$ and
$\varphi'$. By additivity, using affine charts, the previous
construction may be extended to any morphism $f : S \rightarrow S'$
in $\GDef_k$, in order to define
the notion of
$f$-integrability for
a volume form $\alpha$ in
$|\tilde \Omega|_+ (S)$, and also, when
$\alpha$ is $f$-integrable,
the fiber
integral $f^{\rm top}_! (\alpha)$, which belongs to $|\tilde
\Omega|_+ (S')$. When  $S = h [0, 0, 0]$, we shall say
integrable instead of  $f$-integrable, and we shall write
$\int_S \alpha$ for $f^{\rm top}_! (\alpha)$.

In this framework, one may deduce from (A1b) in Theorem \ref{mt}
the following general form of Fubini Theorem for motivic
integration:

\begin{theorem}[Fubini Theorem]\label{if}Let  $f : S \rightarrow S'$
be a morphism in  $\GDef_k$. Assume  $S$ is of  dimension $s$,
$S'$ is of  dimension $s'$, and that the fibers of  $f$ are all of
dimension $s - s'$. A positive volume form  $\alpha$ in
$|\tilde \Omega|_+ (S)$ is integrable if and only if it is
$f$-integrable and $f^{\rm top}_! (\alpha)$ is integrable. When
this holds, then
\begin{equation}
\int_S \alpha = \int_{S'} f^{\rm top}_! (\alpha).
\end{equation}
\end{theorem}

\subsection{Comparison with classical motivic integration}\label{3.1}In the definition
of  $\Def_k$, $\RDef_k$ and $\GDef_k$, instead of considering the
category $F_k$ of all fields containing $k$, one could as well
restrict to the subcategory ${\rm ACF}_k$ of algebraically closed
fields containing $k$ and define categories $\Def_{k, {\rm
ACF}_k}$, etc. In fact, it is a direct consequence of Chevalley's
constructibility theorem that $K_0 (\RDef_{k, {\rm ACF}_k})$ is
nothing else than  the  Grothendieck ring $K_0 ({\rm Var}_k)$
considered in   \cite{arcs}. It follows that there is a canonical
morphism $SK_0 (\RDef_{k}) \rightarrow K_0 ({\rm Var}_k)$ sending
$\LL$ to the class of  $\AA^1_k$, which we shall still denote by
$\LL$. One can extend this morphism to a morphism $\gamma : SK_0
(\RDef_k) \otimes_{\NN [\LL - 1]} A_+ \rightarrow K_0 ({\rm
Var}_k) \otimes_{\ZZ [\LL]} A$. By considering the series
expansion of $(1 - \LL^{-i})^{-1}$, one defines a canonical
morphism $\delta : K_0 ({\rm Var}_k)\otimes_{\ZZ [\LL]} A
\rightarrow \widehat \cM$, with  $\widehat \cM$ the completion of
$K_0 ({\rm Var}_k) [\LL^{-1}]$ considered in  \cite{arcs}.

Let $X$ be an algebraic variety over $k$ of dimension $d$. Set
$\cX^0 := X \otimes_{\Spec k} \Spec k [[t]]$ and $\cX := \cX^0
\otimes_{\Spec k [[t]]} \Spec k ((t))$. Consider a definable
subassignment $W$ of  $h [\cX]$ in the language  $\LPre$, with the
restriction that constants in the valued field sort that appear in
formulas defining  $W$ in affine charts defined over  $k$ belong
to $k$ (and not to  $k ((t))$). We assume  $W (K) \subset \cX (K
[[t]])$ for every  $K$ in $F_k$. With the notation of \cite{arcs},
formulas defining $W$ in affine charts define a semialgebraic
subset of the arc space $\cL (X)$ in the corresponding chart, by
Theorem \ref{quantel} and Chevalley's constructibility theorem. In
this way we assign canonically to $W$ a semialgebraic subset
$\tilde W$ of $\cL (X)$. Similarly, let $\alpha$ be a definable
function on  $W$ taking integral values and satisfying the
additional condition that constants in the valued field sort,
appearing in formulas defining $\alpha$ can only belong to $k$. To
any such function $\alpha$ we may assign a semialgebraic function
$\tilde \alpha$ on  $\tilde W$.

\begin{theorem}\label{compvar}Under the former hypotheses,
$|\omega_0|$ denoting the canonical volume form on  $h [\cX]$, for
every definable function $\alpha$ on $W$ with integral values
satisfying the previous conditions and bounded below, $\11_W
\LL^{-\alpha}|\omega_0|$ is integrable on $h [\cX]$ and
\begin{equation}
(\delta \circ \gamma) \Bigl(\int_{h [\cX]} \11_W \LL^{-\alpha}|\omega_0|\Bigr)
=
\int_{\tilde W} \LL^{- \tilde \alpha} d \mu',
\end{equation}
$\mu'$ denoting the motivic measure considered in   \cite{arcs}.
\end{theorem}

It follows from Theorem \ref{compvar} that, for semialgebraic sets
and functions, the motivic integral constructed in \cite{arcs} in
fact already exists in  $K_0 ({\rm Var}_k) \otimes_{\ZZ [\LL]} A$,
or even in $SK_0 ({\rm Var}_k)\otimes_{\NN [\LL-1]} A_+$, with
$SK_0 ({\rm Var}_k)=SK_0 (\RDef_{k, {\rm ACF}_k})$, the
Grothendieck semiring of varieties over $k$.

\subsection{Comparison with arithmetic motivic integration}\label{4545}
Similarly, instead of ${\rm ACF}_k$, we may also consider the category ${\rm
PFF}_k$ of pseudo-finite fields containing $k$. Let us recall that
a pseudo-finite field is a perfect field  $F$ having a unique extension of
degree  $n$ for every  $n$ in a given algebraic closure and such
that every geometrically irreducible variety over $F$ has a
$F$-rational point. By restriction from $F_k$ to ${\rm PFF}_k$ we can
define categories $\Def_{k, {\rm PFF}_k}$, etc. In particular,
the  Grothendieck ring $K_0 (\RDef_{k, {\rm PFF}_k})$ is nothing
else but what is denoted by $K_0 ({\rm PFF}_k)$ in  \cite{pek} and
\cite{dw}.

In the paper \cite{JAMS},
arithmetic motivic integration was taking its values in a certain completion
$\hat K_0^v ({\rm Mot}_{k, \bar \QQ})_{\QQ}$ of a ring
$K_0^v ({\rm Mot}_{k, \bar \QQ})_{\QQ}$.
Somewhat later it was remarked in
\cite{pek} and  \cite{dw} one can restrict to
the smaller ring $K_0^{\rm mot} ({\rm Var}_k) \otimes \QQ$,
the definition of which we shall now recall.

The field $k$ being of characteristic $0$, there exists, by
\cite{GS} and \cite{GN}, a unique morphism of rings $K_0 ({\rm
Var}_k) \rightarrow K_0 ({\rm CHMot}_k)$ sending the class of a
smooth projective variety $X$ over $k$ to the class of its Chow
motive. Here $K_0 ({\rm CHMot}_k)$ denotes the Grothendieck ring
of the category of Chow motives over $k$ with rational
coefficients. By definition, $K_0^{\rm mot} ({\rm Var}_k) $ is the
image of $K_0 ({\rm Var}_k)$ in $K_0 ({\rm CHMot}_k)$ under this
morphism. [Note that the definition of  $K_0^{\rm mot} ({\rm
Var}_k) $ given in   \cite{pek} is  not clearly equivalent and
should be replaced by the one given above.] In \cite{pek} and
\cite{dw}, the authors have constructed, using results from
\cite{JAMS}, a canonical morphism $\chi_c : K_0 ({\rm PFF}_k )
\rightarrow K_0^{\rm mot} ({\rm Var}_{k}) \otimes \QQ$ as follows:

\begin{theorem}[Denef-Loeser \cite{pek} \cite{dw}]\label{imp}Let $k$
be a field of characteristic zero.
There exists a unique ring morphism
\begin{equation}
\chi_c : {\rm K}_0
({\rm PFF}_k ) \longrightarrow {\rm K}_0^{\rm mot} ({\rm Var}_{k}  )
\otimes {\bf Q}
\end{equation}
satisfying  the following two properties:
\begin{enumerate}\item[(i)] For any formula $\varphi$ which is a conjunction of polynomial equations
over $k$, the element $\chi_c ([\varphi])$ equals the class in
${\rm K}_0^{\rm mot} ({\rm Var}_{k}  ) \otimes {\bf Q}$ of the
variety defined by $\varphi$. \item[(ii)] Let $X$ be a normal
affine irreducible variety over $k$, $Y$ an unramified Galois
cover of $X$, that is, $Y$ is an integral \'etale scheme over $X$
with $Y/G \cong X$, where $G$ is the group of all endomorphisms of
$Y$ over $X$, and $C$ a cyclic subgroup of the Galois group G of
$Y$ over $X$. For such data we denote by $\varphi_{Y,X,C}$ a ring
formula whose interpretation, in any field $K$ containing $k$, is
the set of $K$-rational points on $X$ that lift to a geometric
point on $Y$ with decomposition group $C$ (i.e., the set of points
on $X$ that lift to a $K$-rational point of $Y/C$, but not to any
$K$-rational point of $Y/C'$ with $C'$ a proper subgroup of $C$).
Then
\[
\chi_c ([\varphi_{Y,X,C} ]) = \frac{{\left| C \right|}}{{\left|
{{\rm N}_{G} (C)} \right|}}\chi_c ([\varphi_{Y,Y/C,C}
]),\]
where ${\rm N}_{G} (C)$ is the normalizer of $C$ in $G$.
\end{enumerate}
Moreover, when  $k$ is a number field, for almost all finite
places $\gP$, the number of rational points of $(\chi_c
([\varphi]))$ in the residue field $k(\gP)$ of $k$ at $\gP$ is
equal to the cardinality of $h_{\varphi} (k(\gP))$.
\end{theorem}

The construction of $\chi_c$ has been recently extended
to the relative setting by J. Nicaise \cite{nicaise}.

\subsection{}\label{4646} The arithmetical measure takes its values
in a certain completion $ \hat K_0^{\rm mot} ({\rm Var}_{k})
\otimes \QQ$ of the localisation of $K_0^{\rm mot} ({\rm Var}_{k})
\otimes \QQ$ with respect to the class of the affine line. There
is a canonical morphism $\hat \gamma : SK_0 (\RDef_{k})
\otimes_{\NN [\LL- 1]} A_+ \rightarrow K_0 ({\rm PFF}_k)
\otimes_{\ZZ [\LL]} A$.  Considering the series expansion of $(1 -
\LL^{-i})^{-1}$, the map $\chi_c$ induces a canonical morphism
$\tilde \delta : K_0 ({\rm PFF}_k) \otimes_{\ZZ [\LL]} A
\rightarrow \hat K_0^{\rm mot} ({\rm Var}_{k} )\otimes \QQ$.

Let  $X$ be an algebraic variety over  $k$ of  dimension $d$. Set
$\cX^0 := X \otimes_{\Spec k} \Spec k [[t]]$, $\cX := \cX^0
\otimes_{\Spec k [[t]]} \Spec k ((t))$, and consider a definable
subassignment $W$ of   $h [\cX]$ satisfying the conditions in
\ref{3.1}. Formulas defining $W$ in affine charts allow to define,
in the terminology and with the notation in \cite{JAMS}, a
definable subassignment of $h_{\cL (X)}$ in the corresponding
chart, and we may assign canonically to $W$ a definable
subassignment  $\tilde W$ of  $h_{\cL (X)}$ in the sense of
\cite{JAMS}.

\begin{theorem}\label{arcompvar}Under the previous hypotheses and with the previous notations,
$\11_W |\omega_0|$ is integrable on
$h [\cX]$ and
\begin{equation}
(\tilde \delta \circ \hat \gamma) \Bigl(\int_{h [\cX]} \11_W
|\omega_0|\Bigr) = \nu (\tilde W),
\end{equation}
$\nu$ denoting the arithmetical motivic measure as defined in
\cite{JAMS}.
\end{theorem}

In particular, it  follows from Theorem \ref{arcompvar} that in the
present setting the arithmetical motivic integral constructed  in
\cite{JAMS} already exists in $K_0 ({\rm PFF}_k) \otimes_{\ZZ
[\LL]} A$ (or even in $SK_0 ({\rm PFF}_k) \otimes_{\NN [\LL-1]}
A_+$), without completing further the Grothendieck ring and
without considering Chow motives (and even without inverting
additively all elements of the Grothendieck semiring).

\section{Comparison with $p$-adic integration}\label{padic}

In the next two sections we present new results on specialization
to $p$-adic integration and
Ax-Kochen-Er{\v s}ov Theorems for integrals
with parameters. We plan to give complete details in a future paper.

\subsection{$P$-adic definable sets}\label{17.1}We fix a finite extension
$K$ of $\QQ_p$ together with an uniformizing parameter $\varpi_K$.
We denote by $R_K$ the  valuation ring and by $k_K$ the residue
field, $k_K \simeq \FF_{q(K)}$ for some power $q(K)$ of $p$. Let
$\varphi$ be a formula in the language $\LPre$ with coefficients
in $K$ in the valued field sort and coefficients in $k_K$ in the
residue field sort, with $m$ free variables in the valued field
sort, $n$ free variables in the residue field sort and $r$ free variables
in the value group
sort. The formula $\varphi$ defines a subset $Z_{\varphi}$ of $K^m
\times k_K^n \times \ZZ^r$ (recall that since we have chosen
$\varpi_K$, $K$ is endowed with an angular component mapping). We
call such a subset a $p$-adic definable subset of $K^m \times
k_K^n \times \ZZ^r$. We define morphisms between $p$-adic
definable subsets similarly as before: if $S$ and $S'$ are
$p$-adic definable subsets of $K^m \times k_K^n \times \ZZ^r$ and
$K^{m'} \times k_K^{n'} \times \ZZ^{r'}$ respectively, a morphism
$f : S \rightarrow S'$ will be a function $f : S \rightarrow S'$
whose graph is $p$-adic definable.

\subsection{$P$-adic dimension}\label{17.2}By the work of
Scowcroft and van den Dries \cite{sd}, there is a good dimension
theory for $p$-adic definable subsets of $K^m$. By Theorem 3.4 of
\cite{sd}, a  $p$-adic definable subset $A$ of $K^m$ has dimension
$d$ if and only its Zariski closure has dimension $d$ in the sense
of algebraic geometry. For $S$ a $p$-adic definable subset of $K^m
\times k_K^n \times \ZZ^r$, we define the dimension of $S$ as the
dimension of its image $S'$ under the projection $\pi : S
\rightarrow K^m$.

More generally if $f : S \rightarrow S'$ is a morphism of $p$-adic
definable subsets, one defines the relative dimension of $f$ to be
the maximum of the dimensions of the fibers of $f$.

\subsection{Functions}\label{17.3}
Let $S$ be a $p$-adic definable subset of $K^m \times k_K^n \times
\ZZ^r$. We shall consider the $\QQ$-algebra $\cC_{K} (S)$
generated by functions of the form $\alpha$ and $q^{\alpha}$ with
$\alpha$ a $\ZZ$-valued $p$-adic definable function on $S$. For
$S'\subset S$ a $p$-adic definable subset, we write $\11_{S'}$ for
the characteristic function of $S'$ in $\cC_{K} (S)$.

For $d \geq 0$ an integer, we  denote by $\cC_{K}^{\leq d}  (S)$
the ideal of  $\cC(S)$ generated by all functions $\11_{S'}$ with
$S'$ a $p$-adic definable subset of $S$ of dimension $\leq d$.
Similarly to what we did before, we set
\begin{equation}
C_{K}^d (S) := \cC_{K}^{\leq d}  (S) /\cC_{K}^{\leq d- 1}  (S) \quad
\text{and} \quad C_{K} (S) := \bigoplus_d C_{K}^d (S).
\end{equation}

Also, similarly as before, we have relative
variants of the above definitions.
If $f : Z \rightarrow S$
is a morphism
between
$p$-adic definable subsets,
we define
$\cC_{K}^{\leq d}  (Z \rightarrow S)$,
$C_{K}^d (Z \rightarrow S)$ and $C_{K} (Z \rightarrow S)$
by replacing dimension by relative dimension.

\subsection{$P$-adic measure}\label{17.4} Let $S$ be a $p$-adic definable subset of $K^m$
of dimension $d$. By the construction of \cite{Veys} based on
\cite{serre}, bounded $p$-adic definable subsets $A$ of $S$ have a
canonical $d$-dimensional volume $\mu_{K}^d (A)$ in $\RR$.

Now let $S$ be a $p$-adic definable subset of $K^m\times
k_K^n\times \ZZ^r$ of dimension $d$ and $S'$ its image under the
projection $\pi : S \rightarrow K^m$. We define the measure
$\mu_d$ on $S$ as the measure induced by the product  measur on
$S'\times k_K^n \times \ZZ^r$ of
the
$d$-dimensional volume $\mu_{K}^d$ on the factor $S'$
and  the counting measure on the factor $k_K^n\times
\ZZ^r$. When $S$ is of dimension $<d$ we
declare $\mu_{K}^d$ to be identically zero.

We call $\varphi$ in $\cC_{K} (S)$ integrable on $S$ if $\varphi$
is integrable against $\mu_d$ and we denote the integral by
$\mu_{K}^d(\varphi)$.

One defines ${\rm I} C_{K}^d (S)$ as the abelian subgroup of
$C_{K}^d (S)$ consisting of the classes of integrable functions in
$\cC_{K} (S)$. The measure $\mu_{K}^d$ induces a morphism of
abelian groups $\mu_{K}^d : {\rm I} C_{K}^d (S) \rightarrow \RR$.

More generally if $\varphi=\varphi\11_{S'}$, where $S'$ has
dimension $i \leq d$, we say $\varphi$ is $i$-integrable if its
restriction $\varphi'$ to $S'$ is integrable and we set $\mu_{K}^i
(\varphi) := \mu_{K}^i (\varphi')$. One defines  ${\rm I} C_{K}^i
(S)$ as the abelian subgroup of $C_{K}^i (S)$ of the classes of
$i$-integrable functions in $\cC_{K} (S)$. The measure $\mu_{K}^i$
induces a morphism of abelian groups $\mu_{K}^i : {\rm I} C_{K}^i
(S) \rightarrow \RR$. Finally we set ${\rm I} C_{K} (S) :=
\bigoplus_i {\rm I} C_{K}^i (S)$ and we define $\mu_{K} : {\rm I}
C_{K} (S) \rightarrow \RR$ to be the sum of the morphisms
$\mu_{K}^i$. We call elements of $C_{K} (S)$, resp. ${\rm I} C_{K}
(S)$, constructible Functions, resp. integrable constructible
Functions on $S$.

Also, if $f : S \rightarrow \Lambda$ is a morphism of $p$-adic
definable subsets, we shall say an element $\varphi$ in $C_{K} (S
\rightarrow \Lambda)$ is integrable if the restriction of
$\varphi$ to every fiber of $f$ is an integrable constructible
 Function and we denote by ${\rm I} C_{K} (S  \rightarrow \Lambda)$
the set of such Functions.

We may now reformulate Denef's basic Theorem on $p$-adic integration
(Theorem 1.5 in \cite{D2000}, see also \cite{D85}):

\begin{theorem}[Denef]\label{dbt}Let $f : S \rightarrow \Lambda$
be a morphism of $p$-adic definable subsets. For every integrable
constructible Function $\varphi$ in $C_{K} (S \rightarrow
\Lambda)$, there exists a unique function $\mu_{K, \Lambda}
(\varphi)$ in $\cC (\Lambda)$ such that, for every point $\lambda$
in $\Lambda$,
\begin{equation}
\mu_{K, \Lambda} (\varphi) (\lambda) =
\mu_{K} (\varphi_{|f^{-1} (\lambda)}).
\end{equation}
\end{theorem}

Strictly speaking, this is not the statement that one finds
in \cite{D2000}, but the proof
sketched there extends to our setting.

\subsection{Pushforward}\label{sfr}
It is possible to define, for every morphism $f : S \rightarrow S'$
of
$p$-adic definable subsets,
a natural pushforward morphism
\begin{equation}
f_! : {\rm I} C_{K} (S) \longrightarrow {\rm I} C_{K} (S')
\end{equation}
satisfying similar properties
as in Theorem \ref{mt}.
This may be done along similar lines as what we did
in the motivic case
using Denef's $p$-adic cell decomposition \cite{Dcell} instead of Denef-Pas
cell decomposition. Note however that much  less work is required
in this case, since one already knows what the $p$-adic measure is!
In particular, when $f$ is the projection on the one point definable subset
one recovers the $p$-adic measure $\mu_K$.
Also in the relative setting we have
natural pushforward morphisms
\begin{equation}
f_{!\Lambda} : {\rm I} C_{K} (S \rightarrow \Lambda)
\longrightarrow {\rm I} C_{K} (S' \rightarrow \Lambda),
\end{equation}
for $f : S \rightarrow S'$ over $\Lambda$, and one recovers the
relative $p$-adic measure $\mu_{K, \Lambda}$ when $f$ is the
projection to $\Lambda$.

\subsection{Comparison with $p$-adic integration}\label{nf}
Let $k$ be a number field with ring of integers $\cO$. Let
$\cA_{\cO}$ be the collection of all the $p$-adic completions of
$k$ and of all finite field extensions of $k$. In this section and
in section \ref{red} we let $\LO$ be the language $\LPre (\cO
[[t]])$, that is, the language $\LPre$ with coefficients in $k$
for the residue field sort and coefficients in $\cO [[t]]$ for the
valued field sort, and, all definable subassignments, definable
morphisms, and motivic constructible functions will be with
respect to this language.
To stress the fact that our language is
$\LO$ we use the notation $\Def(\LO)$ for $\Def$, and similarly
for $\cC (S, \LO)$, $\Def_S (\LO)$ and so on.

For $K$ in $\cA_{\cO}$ we write $k_K$ for its residue field with
$q(K)$ elements, $R_K$ for its valuation ring and $\varpi_K$ for a
uniformizer of $R_K$.

Let us choose  for a while,
for every definable subassignment $S$ in $\Def(\LO)$, a
$\LO$-formula $\psi_S$ defining $S$.
We shall write $\tau(S)$ to
denote the datum $(S,\psi_S)$. Similarly, for any element $\varphi$ of
$\cC(S)$, $C(s)$, $IC(S)$, and so on, we choose a finite set $\psi_{\varphi, i}$
of formulas needed to determine $\varphi$ and we
write $\tau(\varphi)$ for $(\varphi,\{\psi_{\varphi,i}\}_i)$.

Let $S$ be a definable subassignment of $h [m, n, r]$ in
$\Def(\LO)$ with $\tau (S) = (S,\psi_S)$.
Let  $K $ be in $\cA_{\cO}$. One may consider $K$ as an $\cO[[t]]$-algebra
via
the morphism
\begin{equation}
\lambda_{\cO,K}:\cO[[t]]\to K:\sum_{i\in\NN}a_it^i\mapsto
\sum_{i\in\NN}a_i\varpi_K^i,
\end{equation}
hence, if one interprets elements $a$ of
$\cO[[t]]$ as $\lambda_{\cO,K}(a)$,
the formula $\psi_S$ defines a $p$-adic definable subset $S_{K,\tau}$ of $K^m
\times k_K^n \times \ZZ^r$.

If now $\tau (S) = (S,\psi_S)$ is replaced by $\tau' (S) =
(S,\psi'_S)$ with $\psi'_S$ another $\LO$-formula defining $S$, it
follows,  from a small variant of Proposition 5.2.1 of \cite{JAMS}
(a result of Ax-Kochen-Er{\v s}ov type that uses ultraproducts and
follows from the Theorem of Denef-Pas), that there exists an
integer $N$ such that $S_{K,\tau}=S_{K,\tau'}$ for every $K$ in
$\cA_{\cO}$ with residue field characteristic $\chara k_K \geq N$.
(Note however that this number $N$ can be arbitrarily large for
different $\tau'$.)

Let us consider the quotient
\begin{equation}
\prod_{K \in \cA_{\cO}}\cC_{K}
(S_{K, \tau}) /\ \sum_N \prod_{\sur{K \in \cA_{\cO}}{\chara k_K<N}} \cC_{K}
(S_{K, \tau}),
\end{equation}
consisting of  families indexed by $K$
of elements of $\cC_{K} (S_{K, \tau})$, two such
families being  identified if for some $N > 0$ they coincide for
$\chara k_K \geq N$. It follows from the above remark that it
is independent of $\tau$ (more precisely all these quotients are canonically isomorphic),
so we may denote it by
\begin{equation}
\prod{}^{'} \cC_{K} (S_{K}).
\end{equation}
 One defines similarly $\prod{}^{'} C_{K}
(S_{K})$, $\prod{}^{'} {\rm I}C_{K} (S_{K})$, etc.

Now take $W$ in $\RDef_S (\LO)$. It defines a $p$-adic definable
subset $W_{K, \tau}$ of $S_{K, \tau} \times (k_K)^{\ell}$, for some $\ell$,
for every  $K$ in $ \cA_{\cO}$. We
may now consider the function $\psi_{W, K, \tau}$ on $S_{K, \tau}$ assigning
to a
point $x$ the number of points mapping to it in $W_{K, \tau}$,
that is, $\psi_{W, K, \tau} (x) = {\rm card} (W_{K, \tau} \cap (\{x\} \times
k_K^{\ell}))$.
Similarly as before, if we take another
function $\tau'$, we have $\psi_{W, K, \tau} = \psi_{W, K, \tau'}$ for every $K$ in $ \cA_{\cO}$ with
residue field characteristic $\chara k_K \geq N$,
hence we get in this way an arrow $\RDef_S (\LO) \rightarrow \prod{}^{'}
\cC_{K} (S_{K})$ which factorizes through a ring morphism $K_0
(\RDef_S (\LO)) \rightarrow \prod{}^{'} \cC_{K} (S_{K})$. If we
send $\LL$ to $q(K)$, one can extend uniquely this morphism to a
ring morphism
\begin{equation}
\Gamma : \cC (S, \LO) \longrightarrow \prod{}^{'}\cC_{K} (S_{K}).
\end{equation}

Since $\Gamma$ preserves the (relative) dimension of support on
those factors $K$ with $\chara k_K$ big enough, $\Gamma$ induces the
morphisms
\begin{equation}
\Gamma : C (S, \LO) \longrightarrow \prod{}^{'} C_{K} (S_{K})
\end{equation}
and
\begin{equation}
\Gamma : C (S \rightarrow \Lambda,  \LO) \longrightarrow
\prod{}^{'} C_{K} (S_{K} \rightarrow \Lambda_{K}),
\end{equation}
for $S \rightarrow \Lambda$ a morphism in $\Def_K (\LO)$.

The
following comparison Theorem says that
the morphism $\Gamma$ commutes with pushforward. In more concrete terms,
given an integrable function $\varphi$
in
$C (S \rightarrow \Lambda, \LO)$, for almost all $p$,
its specialization $\varphi_K$
to any finite extension $K$
of $\QQ_p$ in $\cA_{\cO}$ is integrable,
and the specialization of the pushforward of $\varphi$ is equal to
the pushforward of $\varphi_K$.

\begin{theorem}\label{compres0}
Let $\Lambda$ be in $\Def_K (\LO)$ and let $f : S \rightarrow S'$
be a morphism in $\Def_{\Lambda} (\LO)$. The morphism
\begin{equation}\Gamma : C (S \rightarrow \Lambda, \LO) \rightarrow
\prod{}^{'} C_{K} (S_{K} \rightarrow \Lambda_{K})
\end{equation}
induces a morphism
\begin{equation}
\Gamma: {\rm I}C (S \rightarrow \Lambda, \LO)  \to \prod{}^{'}
{\rm I}C_{K} (S_{K} \rightarrow \Lambda_{K})
\end{equation}
(and similarly
for $S'$), and the following diagram is commutative:
\begin{equation*}\label{777}\xymatrix{
{\rm I}C (S \rightarrow \Lambda, \LO) \ar[r]^>>>>{\Gamma}
\ar[d]_{f_{!\Lambda}} & {\prod{}^{'}} {\rm I}C_{K}
(S_{K}\rightarrow \Lambda_{K})
\ar[d]^{\prod{}^{'} f_{K, \Lambda_{K}!}}\\
 {\rm I}C (S' \rightarrow \Lambda, \LO)
\ar[r]^>>>>{\Gamma}& {\prod{}^{'} }{\rm I}C_{K} (S'_{K}\rightarrow
\Lambda_{K}), }
\end{equation*}
with $f_{K} : S_{K} \rightarrow S_{K}'$ the morphism induced by
$f$ and where the map $\prod{}^{'} f_{K, \Lambda_{K}!} $ is
induced by the maps
 $f_{K, \Lambda_{K}!}:{\rm I}C_{K}
(S_{K}\rightarrow \Lambda_{K})\to {\rm I}C_{K} (S'_{K}\rightarrow
\Lambda_{K}).$
\end{theorem}
\begin{proof}[Sketch of proof] The image of $\varphi$ in ${\rm I}C (S \rightarrow \Lambda,
\LO)$ under $f_{!\Lambda}$ can be calculated by taking an
appropriate cell decomposition of the occurring sets, adapted to
the occurring functions (as in \cite{cl} and inductively applied
to all valued field variables). Such calculation is independent of
the choice of cell decomposition by the unicity statement of
Theorem \ref{mt}. By the Ax-Kochen-Er{\v s}ov principle for the
language $\LO$ implied by Theorem \ref{quantel}, this cell
decomposition determines, for $K$ in $\cA_\cO$ with $\chara k_K$
sufficiently large, a cell decomposition \`{a} la Denef (in the
formulation of Lemma 4 of \cite{C}) of the $K$-component of these
sets, adapted to the $K$-component of the functions occuring here,
where thus the same calculation can be pursued. That this
calculation is actually the same follows from the fact that
$p$-adic integration satisfies properties analogue to the  axioms
of Theorem \ref{mt}.
\end{proof}

In particular, we have the following statement, which says that,
given an integrable function $\varphi$ in $C (S \rightarrow
\Lambda, \LO)$, for almost all $p$, its specialization $\varphi_F$
to any finite extension $F$ of $\QQ_p$ in $\cA_{\cO}$ is
integrable, and the specialization of the motivic integral $\mu
(\varphi)$ is equal to the $p$-adic integral of $\varphi_F$:

\begin{theorem}\label{compres1}Let
$f : S \rightarrow \Lambda$ be a morphism in $\Def_K (\LO)$.
The following diagram is commutative:
\begin{equation*}\label{44444}\xymatrix{
{\rm I}C (S \rightarrow \Lambda, \LO) \ar[r]^>>>>{\Gamma}
\ar[d]_{\mu_{\Lambda}} & {\prod{}^{'} } {\rm I}C_{K} (S_{K}
\rightarrow \Lambda_{K})
\ar[d]^{\prod{}^{'}\mu_{K, \Lambda_{K}}}\\
 \cC (\Lambda, \LO)
\ar[r]^<<<<<<<<{\Gamma}& {\prod{}^{'}} \cC_{K} (\Lambda_{K}). }
\end{equation*}
\end{theorem}

\section{Reduction mod $p$ and
a motivic Ax-Kochen-Er{\v s}ov Theorem for integrals
with parameters}\label{18}

\subsection{Integration over $ \FF_q ((t))$}Consider
now the field $K = \FF_q ((t))$ with valuation ring $R_K$ and
residue field $k_K = \FF_{q}$ with $q=q(K)$ a prime power. One may
define $\FF_q ((t))$-definable sets similarly as in \ref{17.1}.
Little is known about the structure of these $\FF_q
((t))$-definable sets, but, for any subset $A$ of $K^m$, not
necessarly definable, we may still define the dimension of $A$ as
the dimension of its Zariski closure. Similarly as in \ref{17.2},
one extends that  definition to any subset $A$ of
 $K^m \times k_K^n \times \ZZ^r$ and define the relative dimension
 of a mapping $f : A \rightarrow \Lambda$, with
 $\Lambda$ any subset of $K^{m'} \times k_K^{n'} \times \ZZ^{r'}$.
 When $A$ is $\FF_q ((t))$-definable, one can define
 a $\QQ$-algebra $\cC_K (A)$ as in
 \ref{17.3}, but since no analogue of Theorem \ref{dbt}
 is known in this setting, we shall
 consider, for $A$ any subset of
 $K^m \times k_K^n \times \ZZ^r$,
 the $\QQ$-algebra $\cF_K (A)$ of all functions
 $A \rightarrow \QQ$. For $d \geq 0$ an integer, we denote
 by $\cF_K^{\leq d}(A)$ the ideal of functions with support of dimension $\leq d$.
 We set $F^d_K (A):= \cF_K^{\leq d}(A) / \cF_K^{\leq d- 1}(A)$
 and
 $F_K (A) := \oplus_d F_K^d (A)$.
 One defines similarly relative variants
 $\cF^{\leq d}_K (A\rightarrow \Lambda)$, $F^d_K (A\rightarrow \Lambda)$ and
 $F_K (A\rightarrow \Lambda)$,
 for $f : A \rightarrow A'$ as above.

Let $A$ be a subset of $K^m$ with Zariski closure $\bar A$ of
dimension $d$. We consider the canonical $d$-dimensional measure
$\mu_K^d$ on $\bar A (K)$ as in \cite{oesterle}. We say a function
$\varphi$ in $\cF_K (A)$ is integrable if it is measurable and
integrable with respect to the measure $\mu_K^d$. Now we may
proceed as in \ref{17.4} to define, for $A$ a subset of $K^m
\times k_K^n \times \ZZ^r$, ${\rm I} F_K (A)$ and $\mu_K : {\rm I}
F_K (A) \rightarrow \RR$. Also, if $f : A \rightarrow \Lambda$ is
a mapping as before, one defines ${\rm I} F_K (A \rightarrow
\Lambda)$ as Functions whose restrictions to all fibers lie in
${\rm I} F_K$. We denote by $\mu_{K, \Lambda}$ the unique mapping
${\rm I} F_K (A \rightarrow \Lambda) \rightarrow \cF (\Lambda)$
such that, for every $\varphi$ in ${\rm I} F_K (A \rightarrow
\Lambda)$ and every point $\lambda$ in $\Lambda$, $ \mu_{K,
\Lambda} (\varphi) (\lambda) = \mu_{K} (\varphi_{|f^{-1}
(\lambda)})$.

\subsection{Reduction mod $p$}\label{red} We go back to the notation of \ref{nf}.
In particular, $k$ denotes a number field with ring of integers
$\cO$, $\cA_{\cO}$ denotes the set of all $p$-adic completions of
$k$ and of all the finite field extensions of $k$, and $\LO$
stands for the language $\LPre (\cO [[t]])$. We also use the map
$\tau$ as defined in section \ref{nf}.

Let $\cB_{\cO}$ be the set of all local fields  over $\cO$ of
positive characteristic. As for $\cA_{\cO}$, we use for every  $K$
in
$\cB_{\cO}$ the notation $k_K$ for its residue field with $q(K)$
elements, $R_K$ for its valuation ring and $\varpi_K$ for a
uniformizer of $R_K$.

Let $S$ be a definable subassignment of $h [m, n, r]$ in
$\Def(\LO)$ and let $\tau(S)$ be $(S,\psi_S)$ with $\psi_S$ a
$\LO$-formula. Similarly as for $\cA_{\cO}$, since every
$K$ in $\cB_{\cO}$ is an $\cO[[t]]$-algebra under the morphism
\begin{equation}
\lambda_{\cO,K}:\cO[[t]]\to K:\sum_{i\in\NN} a_it^i \mapsto
\sum_{i\in\NN}a_i\varpi_K^i,
\end{equation}
interpreting any element $a$ of $\cO[[t]]$ as
$\lambda_{\cO,K}(a)$, $\psi_S$ defines a $K$-definable subset
$S_{K,\tau}$ of $K^m \times k_K^n \times \ZZ^r$. Again by a small
variant of Proposition 5.2.1 of \cite{JAMS}, for any other $\tau'$
we have for every  $K$ in $\cB_{\cO}$ with $\chara k_K$ big enough
that $S_{K,\tau}=S_{K,\tau'}$, hence, may define, similarly as in
\ref{nf},
\begin{equation}
\prod{}^{'} \cF_{K} (S_{K}).
\end{equation}
to be
the quotient
\begin{equation}
\prod_{K \in \cB_{\cO}}\cF_{K}
(S_{K, \tau}) /\ \sum_N \prod_{\sur{K \in \cB_{\cO}}{\chara k_K<N}} \cF_{K}
(S_{K, \tau}),
\end{equation} and  similarly  for $\prod{}^{'} F_{K} (S_{K})$, $\prod{}^{'}
{\rm I}F_{K} (S_{K})$, etc.

Similarly as in \ref{nf},
one may define  ring morphisms
\begin{equation}
\hat \Gamma : \cC (S, \LO) \longrightarrow \prod{}^{'} \cF_{K}
(S_{K}),
\end{equation}
\begin{equation}
\hat \Gamma : C (S, \LO) \longrightarrow \prod{}^{'} F_{K}
(S_{K})
\end{equation}
and
\begin{equation}
\hat \Gamma : C (S \rightarrow \Lambda, \LO) \longrightarrow
\prod{}^{'} F_{K} (S_{K} \rightarrow \Lambda_{K}),
\end{equation}
for $S \rightarrow \Lambda$ a morphism in $\Def_K (\LO)$.

The following statement is a companion to Theorem \ref{compres1}
and has an essentially similar proof.

\begin{theorem}\label{compres3}Let
$f : S \rightarrow \Lambda$ be a morphism in $\Def_K (\LO)$. The
morphism
\begin{equation}
\hat \Gamma : C (S \rightarrow \Lambda, \LO)
\rightarrow \prod{}^{'} F_{K} (S_{K} \rightarrow \Lambda_{K})
\end{equation}
induces a morphism
\begin{equation}
\hat \Gamma : {\rm I}C (S \rightarrow \Lambda, \LO) \to
\prod{}^{'} {\rm I} F_{K} (S_{K} \rightarrow \Lambda_{K})
\end{equation}
and the following diagram is commutative:
\begin{equation*}\label{9999}\xymatrix{
{\rm I}C (S \rightarrow \Lambda, \LO) \ar[r]^>>>>{\hat \Gamma}
\ar[d]_{\mu_{\Lambda}} & {\prod{}^{'} } {\rm I}F_{K} (S_{K}
\rightarrow \Lambda_{K} )
\ar[d]^{\prod{}^{'}\mu_{K, \Lambda_{K} }}\\
\cC (\Lambda, \LO)
\ar[r]^<<<<<<<<{\hat \Gamma}& {\prod{}^{'}} \cF_{K}
(\Lambda_{K}). }
\end{equation*}
\end{theorem}

\subsection{Ax-Kochen-Er{\v s}ov Theorems for motivic
integrals}

We keep the notation of \ref{nf} and \ref{red}. Let $S$, resp.
$\Lambda$, be definable subassignments of $h [m, n, r]$, resp. $h
[m', n', r']$, in the language $\LO$ and consider a definable (in
the language $\LO$) morphism $f : S \rightarrow \Lambda$. Since we
are interested in integrals along the fibers of $f$, there is no
restriction in assuming, and we shall do so, that $\Lambda = h
[m', n', r']$. We set $\Lambda (\cO) := \cO [[t]]^{m'} \times
k^{n'} \times \ZZ^{r'}$.

A first attempt to get Ax-Kochen-Er{\v s}ov Theorems for motivic
integrals is by comparing values. This is achieved as follows. To
every point $\lambda$ in $\Lambda (\cO)$ we may assign, for all
$K$ in $\cA_{\cO}\cup\cB_{\cO}$, a point $\lambda_{K}$ in
$(R_K)^{m'}\times k_K^{n'}\times \ZZ^{r'}$, by using the maps
$\lambda_{\cO,K}$ on the $\cO [[t]]^{m'}$-factor and reduction
modulo $\chara k_K$ for the $k^{n'}$-factor.

Let  $\varphi$ be in $C (S \rightarrow \Lambda, \LO)$. With a
slight abuse of notation, we shall write $\varphi_{K}$ for the
component at $K$ in $C_{K} (S_{K} \rightarrow \Lambda_{K})$ of
$\Gamma (\varphi)$, resp.~of $\hat \Gamma (\varphi)$, for $K$ with
$\chara k_K$ big, in $\cA_{\cO}$, resp.~in $\cB_{\cO}$. We shall
use similar notations for $\varphi$ in $\cC (S, \LO) $.

Over the final subassignment $h [0, 0, 0]$ the morphisms
$\Gamma$ and $\hat \Gamma$ have quite simple descriptions.
Indeed, the morphism $\hat\gamma$ of \ref{4646} induces a ring morphism
\begin{equation}
\gamma':\cC ( h_{\Spec k}, \LO)\longrightarrow  K_0 ({\rm PFF}_K)
\otimes_{\ZZ [\LL]} A.
\end{equation}
One the  other hand,
note that $\cC_{K} ({\rm
point}) \simeq \cF_{K'} ({\rm point}) \simeq \QQ$ for $K$ in $\cA_{\cO}$ and $K'$ in
$\cB_{\cO}$. The morphism $\chi_c : K_0 ({\rm PFF}_k ) \rightarrow
K_0^{\rm mot} ({\rm Var}_{k}) \otimes \QQ$ from \ref{4545} induces
a ring morphism
\begin{equation}\delta':K_0 ({\rm PFF}_K) \otimes_{\ZZ [\LL]} A
\to K_0^{\rm mot} ({\rm Var}_{k}) \otimes \QQ \otimes_{\ZZ [\LL]}
A.
\end{equation}
Note that, for  $K$ in $\cA_{\cO}$, resp.~in $\cB_{\cO}$, with
$\chara K$ big enough, the $K$-component of $\Gamma (\alpha)$,
resp.~of $\hat\Gamma(\alpha)$, is equal to the trace of the
Frobenius at $k_K$ acting on an \'etale realisation of
$(\delta'\circ \gamma')(\alpha)$, for $\alpha$ in $\cC ( h_{\Spec
k}, \LO)$. In particular, one deduces the following statement:

\begin{lem}\label{simplelemma}
Let $\psi$ be a function in $\cC (\Lambda, \LO)$.
Then, for every $\lambda$ in $\Lambda (\cO)$,
there exists an integer $N$ such that, for every
$K_1$ in $\cA_{\cO}$, $K_2$ in $\cB_{\cO}$ with $k_{K_1}\simeq k_{K_2}$
and $\chara k_{K_1}>N$,
\begin{equation}
\psi_{K_1}(\lambda_{K_1}) = \psi_{K_2}(\lambda_{K_2}),
\end{equation}
which also is equal to $(i_{\lambda}^* (\psi))_{K_1}$ and to
$(i_{\lambda}^*(\psi))_{K_2}$.
\end{lem}

From Lemma \ref{simplelemma}, Theorem \ref{compres1}
and Theorem \ref{compres3} one deduces immediatly:

\begin{theorem}\label{axk}Let $f : S \rightarrow \Lambda$ be as above.
Let $\varphi$ be a Function in ${\rm I}C (S \rightarrow \Lambda,
\LO)$. Then, for every $\lambda$ in $\Lambda (\cO)$,
there exists an integer $N$ such that for all
$K_1$ in $\cA_{\cO}$, $K_2$ in $\cB_{\cO}$ with $k_{K_1}\simeq k_{K_2}$
and $\chara k_{K_1}>N$
\begin{equation}
\mu_{K_1} (\varphi_{K_1 \vert f^{-1}_{K_1} (\lambda_{K_1})}) =
\mu_{K_2} (\varphi_{K_2
 \vert f_{K_2}^{-1} (\lambda_{ K_2})}),
\end{equation}
which also equals $(\mu_{\Lambda}(\varphi))_{K_1}(\lambda_{K_1})$
and  $(\mu_{\Lambda}(\varphi))_{K_2}(\lambda_{K_2})$.

\end{theorem}

Note that Theorem \ref{thdl} is a corollary of Theorem \ref{axk}
when $(m', n', r') = (0, 0, 0)$. In fact, Theorem \ref{axk} is not
really satisfactory when $(m', n', r') \not= (0, 0, 0)$, since it
is not uniform with respect to $\lambda$. The following example
shows that this unavoidable: take $k = \QQ$, $S = \Lambda = h [1,
0, 0]$, $f$ the identity and $\varphi = \11_{S \setminus \{0\}}$
in ${\rm I} C (S \rightarrow \Lambda) = \cC (S)$. Take $K_1$ in
$\cA_{\cO}$ and $K_2$ in $\cB_{\cO}$.
We have $\varphi_{K_1}(\lambda_{K_1}) =
\varphi_{K_2}(\lambda_{K_2})$ for $\lambda \not=0$ in $\ZZ$ only
if the characteristic of $K_2$ does not divide $\lambda$.

Hence, instead of comparing values of integrals depending on parameters,
we better compare the integrals as functions, which is done as follows:

\begin{theorem}\label{strongaxk}Let $f : S \rightarrow \Lambda$ be as above.
Let $\varphi$ be a Function in ${\rm I}C (S \rightarrow \Lambda,
\LO)$. Then,
there exists an integer $N$ such that for all
$K_1$ in $\cA_{\cO}$, $K_2$ in $ \cB_{\cO}$ with $k_{K_1}\simeq k_{K_2}$
and $\chara k_{K_1}>N$
\begin{equation}
\mu_{K_1, \Lambda_{K_1}} (\varphi_{K_1}) = 0 \quad \text{if and
only if } \quad \mu_{K_2, \Lambda_{K_2}} (\varphi_{K_2}) = 0.
\end{equation}
\end{theorem}
\begin{proof}
Follows directly from Theorem \ref{compres1}, Theorem
\ref{compres3}, and Theorem \ref{strong}.
\end{proof}

\begin{theorem}\label{strong}Let $\psi$ be in
$\cC (\Lambda, \LO)$.
Then,
there exists an integer $N$ such that for all
$K_1$ in $\cA_{\cO}$, $K_2$ in $\cB_{\cO}$ with $k_{K_1}\simeq k_{K_2}$
and $\chara k_{K_1}>N$
\begin{equation}
\psi_{K_1} = 0
\quad
\text{if and only if }
\quad
 \psi_{K_2} = 0.
\end{equation}
\end{theorem}

\begin{remark}Thanks to  results of Cunningham and Hales \cite{CH},
Theorem \ref{strongaxk} applies to the orbital integrals occuring in
the Fundamental Lemma.
Hence, it follows from Theorem \ref{strongaxk} that the Fundamental
Lemma holds over function fields of large characteristic if and only if
it holds for $p$-adic fields of large characteristic.
(Note that the Fundamental Lemma is about the equality of two integrals,
or, which amounts to the same, their difference to be zero.)
In the special situation of the Fundamental
Lemma, a more precise comparison result
has been proved by Waldspurger \cite{wal} by
representation theoretic techniques.
Let us recall that
the Fundamental Lemma for unitary groups has been proved recently by
Laumon and Ng\^o
\cite{ln} for functions fields.
\end{remark}

\bibliographystyle{amsplain}

\begin{thebibliography}{SGA}


\bibitem{AK1}J. Ax and S. Kochen,
\textit{Diophantine problems over local fields. I},
Amer. J. Math. \textbf{87} (1965), 605--630.
\bibitem{AK2}J. Ax and S. Kochen,
\textit{Diophantine problems over local fields. II.: A complete set
of axioms for $p$-adic number theory},
Amer. J. Math. \textbf{87} (1965), 631--648.
\bibitem{AK3}J. Ax and S. Kochen,
\textit{Diophantine problems over local fields. III. Decidable fields},
Ann. of Math. \textbf{83} (1966),  437--456.

\bibitem{C} R. Cluckers: \textit{Classification of semi-algebraic p-adic sets up to semi-algebraic bijection},
 J. Reine Angew. Math. \textbf{540}
(2000) 105 -- 114, math.LO/0311434.


\bibitem{cr1}{R. Cluckers,  F. Loeser},
\textit{Fonctions constructibles et int\'egration motivique I},
math.AG/0403349,
to appear in C. R. Acad. Sci. Paris S\'er. I Math.



\bibitem{cr2}{R. Cluckers,  F. Loeser},
\textit{Fonctions constructibles et int\'egration motivique II},
math.AG/0403350,
to appear in C. R. Acad. Sci. Paris S\'er. I Math.



\bibitem{cl}{R. Cluckers,  F. Loeser},
\textit{Constructible motivic functions and motivic integration},
in preparation.




\bibitem{cohen}P. J. Cohen,
\textit{Decision procedures for real and {$p$}-adic fields},
Comm. Pure Appl. Math.
\textbf{22} (1969),
131--151.



\bibitem{CH}
C. Cunningham, T. Hales,
\textit{Good orbital integrals},
math.RT/0311353.



\bibitem{delon} F. Delon,
\textit{Some $p$-adic model theory},
European women in mathematics (Trieste, 1997), 63--76,
Hindawi Publ. Corp., Stony Brook, NY, 1999.



\bibitem{D85}
J. Denef,
\textit{On the evaluation of certain $p$-adic integrals},
S\'eminaire de th\'eorie des nombres, Paris
1983--84, 25--47, Progr. Math., \textbf{59}, Birkh\"auser
Boston, Boston, MA, 1985


\bibitem{Dcell}
J. Denef,
\textit{$p$-adic semi-algebraic sets and cell decomposition},
J. Reine Angew. Math. \textbf{369} (1986), 154--166.



\bibitem{D2000}
J. Denef,
\textit{Arithmetic and geometric applications of quantifier elimination for valued fields}, Model
theory, algebra, and geometry, 173--198, Math. Sci. Res. Inst. Publ., 39, Cambridge Univ. Press, Cambridge, 2000.




\bibitem{arcs}{J. Denef, F. Loeser},
\textit{Germs of arcs on singular algebraic varieties
and motivic integration},
Invent. Math.
\textbf{135}
(1999),
201--232.

\bibitem{JAMS}
{J. Denef,  F. Loeser},
\textit{Definable sets, motives and $p$-adic integrals},
J. Amer. Math. Soc.,
\textbf{14} (2001), 429--469.




\bibitem{pek}
{J. Denef,  F. Loeser},
\textit{Motivic integration and the {G}rothendieck group of
pseudo-finite fields},
Proceedings of the International Congress of Mathematicians,
Vol. II (Beijing, 2002),
\textbf{2002}, 13--23,
Higher Ed. Press, Beijing.






\bibitem{dw}{J. Denef,  F. Loeser},
\textit{On some rational generating series
occuring in arithmetic geometry}, math.NT/0212202.




\bibitem{vdDries}
{L. van den Dries}, \textit{Dimension of definable sets, algebraic
boundedness and {H}enselian fields}, Ann. Pure Appl. Logic,
\textbf{45} (1989), 189--209.




\bibitem{GS}
{H. Gillet, C. Soul{\'e}},
\textit{Descent, motives and $K$-theory},
J. Reine Angew. Math.
\textbf{478}
(1996),
127--176.

\bibitem{GN}
{F. Guill\'{e}n, V. Navarro Aznar},
\textit{Un crit\`{e}re d'extension d'un foncteur d\'{e}fini sur les
sch\'{e}mas lisses},
Inst. Hautes {\'E}tudes Sci. Publ. Math.
\textbf{95}
(2002), 1--91.


\bibitem{K}
{M. Kontsevich}, Lecture at   Orsay, December
7, 1995.



\bibitem{Ersov}J. Er{\v s}ov,
\textit{On the elementary theory of maximal normed fields},
Dokl. Akad. Nauk SSSR \textbf{165} (1965), 21--23.



\bibitem{ln}
G. Laumon, B.~C. Ng\^o,
\textit{Le lemme fondamental pour les groupes unitaires},
math.AG/0404454.


\bibitem{nicaise}
J. Nicaise,
\textit{Relative motives and the theory of pseudo-finite fields},
math.AG/0403160.

\bibitem{oesterle}
J. Oesterl{\'e}, \textit{R{\'e}duction modulo $p^{n}$ des
sous-ensembles analytiques ferm{\'e}s de ${\ZZ}^{N}_{p}$}, Invent.
Math., \textbf{66} (1982), 325--341.



\bibitem{Pas}
{J. Pas},
\textit{Uniform $p$-adic cell decomposition and local zeta functions},
J. Reine Angew. Math.,
\textbf{399}
(1989),
137--172.



\bibitem{sd}
P. Scowcroft, L. van den Dries,
\textit{On the structure of semialgebraic sets over $p$-adic fields},
J. Symbolic Logic, \textbf{53} (1988), 1138--1164.




\bibitem{serre}
{J.-P. Serre}, \textit{Quelques applications du th\'eor\`eme de densit\'e
de Chebotarev},
Inst. Hautes {\'E}tudes Sci. Publ. Math., \textbf{54} (1981), 323--401.




\bibitem{terj}G. Terjanian,
\textit{Un contre-exemple \`a une conjecture d'Artin},
C. R. Acad. Sci. Paris \textbf{262}  1966 A612.

\bibitem{Veys} W. Veys, \textit{Reduction modulo $p^{n}$ of $p$-adic subanalytic
sets}, Math. Proc. Cambridge Philos. Soc., \textbf{112} (1992),
483--486.



\bibitem{wal}
J.-L. Waldspurger,
\textit{Endoscopie et changement de caract\'{e}ristiques},
preprint, (2004).







\end{thebibliography}

\end{document}